\documentclass[12pt]{article}

\newbox\squ  
\setbox\squ=\hbox{\vrule width.3pt
            \vbox{\hrule height.3pt width.4em\kern1ex\hrule height.3pt}%
            \vrule width.3pt}
\def\endproof{%
  \ifmmode\eqno\copy\squ\smallskip\else{\unskip\nobreak\hfil%
    \penalty50\hskip2em\hbox{}\nobreak\hfil\copy\squ
    \parfillskip=0pt \finalhyphendemerits=0\penalty-100\smallskip}
  \fi}

\usepackage{amsmath}
\usepackage{amssymb}

\newcommand{\scl}{\scriptstyle}
\newcommand{\non}{\nonumber}
\newcommand{\wt}{\widetilde}
\newcommand{\wh}{\widehat}
\newcommand{\ot}{\otimes}
\newcommand{\ve}{\varepsilon}
\newcommand{\lan}{\langle\ts}
\newcommand{\ran}{\ts\rangle}
\newcommand{\ts}{\,}
\newcommand{\Tau}{ {\mathcal T}}
\newcommand{\U}{ {\rm U}}
\newcommand{\Y}{ {\rm Y}}
\newcommand{\CC}{\mathbb{C}}
\newcommand{\ZZ}{\mathbb{Z}}
\newcommand{\A}{\mathcal{A}}
\newcommand{\Cond}{ {\rm Cond}\hspace{1pt}}
\newcommand{\Sym}{\mathfrak S}
\newcommand{\sgn}{ {\rm sgn}\ts}

\newcommand{\gl}{\mathfrak{gl}}
\newcommand{\la}{\lambda}
\newcommand{\La}{\Lambda}
\newcommand{\xiL}{\xi^{}_{\La}}
\newcommand{\xiLo}{\xi^{}_{\La^0}}
\newcommand{\etaL}{\eta^{}_{\La}}
\newcommand{\etaLo}{\eta^{}_{\La^0}}
\newcommand{\al}{\alpha}
\newcommand{\be}{\beta}
\newcommand{\Proof}{\noindent{\it Proof.}\ \ }

\newtheorem{thm}{Theorem}[section]
\newtheorem{prop}[thm]{Proposition}
\newtheorem{lem}[thm]{Lemma}
\newtheorem{cor}[thm]{Corollary}
\newtheorem{defin}[thm]{Definition}
\newtheorem{conj}[thm]{Conjecture}

\newcommand{\bth}{\begin{thm}}
\renewcommand{\eth}{\end{thm}}
\newcommand{\bpr}{\begin{prop}}
\newcommand{\epr}{\end{prop}}
\newcommand{\ble}{\begin{lem}}
\newcommand{\ele}{\end{lem}}
\newcommand{\bco}{\begin{cor}}
\newcommand{\eco}{\end{cor}}
\newcommand{\bconj}{\begin{conj}}
\newcommand{\econj}{\end{conj}}
\newcommand{\bde}{\begin{defin}}
\newcommand{\ede}{\end{defin}}
\newcommand{\beq}{\begin{equation}}
\def\beql#1{\begin{equation}\label{#1}}
\newcommand{\bal}{\begin{aligned}}
\newcommand{\eal}{\end{aligned}}

\begin{document}

\title{\Large\bf Irreducibility criterion for 
tensor products of Yangian evaluation modules}
\author{{\sc A. I. Molev}\\[15pt]
School of Mathematics and Statistics\\
University of Sydney,
NSW 2006, Australia\\
{\tt
alexm\hspace{0.09em}@\hspace{0.1em}maths.usyd.edu.au
}
}

\date{}
\maketitle

\vspace{2 cm}

\begin{abstract}
The evaluation homomorphisms from the Yangian $\Y(\gl_n)$
to the universal enveloping algebra $\U(\gl_n)$ allow one
to regard the irreducible finite-dimensional 
representations of $\gl_n$ as Yangian modules.
We give necessary and sufficient conditions for
irreducibility of tensor products of 
such evaluation modules.
\end{abstract}


\vspace{2 cm}
	


\newpage

\section{Introduction}
\setcounter{equation}{0}

The Yangian $\Y(\gl_n)$ for the general linear Lie algebra $\gl_n$
is a deformation of 
the universal enveloping algebra $\U(\gl_n[x])$ 
in the class of Hopf algebras;
see Drinfeld~\cite{d:ha}.
A theorem of his~\cite{d:nr} (see also Tarasov~\cite{t:im}) provides
a complete discription of finite-dimensional irreducible
representations of $\Y(\gl_n)$ in terms of their highest weights.
Recently, Arakawa~\cite{a:df}
has found a character
formula for each of these representations with the use of the Kazhdan--Lusztig
polynomials; see also Vasserot~\cite{v:aq} for the case of quantum affine
algebras. Nazarov and Tarasov~\cite{nt:ry} (see also Cherednik~\cite{c:ni})
have given an explicit construction for a class of the so-called tame
representations. However, the structure of the general 
finite-dimensional irreducible $\Y(\gl_n)$-module (with $n\geq 3$)
still remains unknown. In this paper we establish 
an irreducibility criterion for
tensor products of the Yangian evaluation modules 
which thus solves this problem
for a wide class of representations of $\Y(\gl_n)$.

Let $\la=(\la_1,\dots,\la_n)$ be an $n$-tuple of complex numbers such that
$\la_i-\la_{i+1}$ is a non-negative
integer for each $i$. Denote by $L(\la)$ 
the irreducible finite-dimensional
representation of the Lie algebra $\gl_n$ with the highest weight $\la$.
For each $a\in\CC$ there is an evaluation homomorphism
$\varphi_a$ from the Yangian $\Y(\gl_n)$ to the universal enveloping algebra
$\U(\gl_n)$: see Section~\ref{sec:pre} for the definitions. 
Using $\varphi_a$ we make
$L(\la)$ into a Yangian module and denote it by $L_a(\la)$.
We keep the notation $L(\la)$ for 
the evaluation module $L_a(\la)$ with $a=0$.
The Hopf algebra structure on $\Y(\gl_n)$ allows one to regard
tensor products	of the type
\beql{tenev}
L_{a_1}(\la^{(1)})\ot L_{a_2}(\la^{(2)})\ot\cdots\ot L_{a_k}(\la^{(k)})
\end{equation}
as Yangian modules. Our main result is a criterion of irreducibility
of these modules: see Theorems~\ref{thm:main1} and 	\ref{thm:main2} below.
To formulate the result we first note that
the problem can be reduced to the particular case
where all the parameters $a_i$ in \eqref{tenev} are equal to zero.
This is done by using the composition of the module \eqref{tenev}
with an appropriate automorphism of the Yangian: 
see Proposition~\ref{prop:redu} below.

We give first an irreducibility criterion for the tensor product
$L(\la)\ot L(\mu)$ of two evaluation modules.
It is well-known (see e.g. \cite{kr:yb} and
Theorem~\ref{thm:sufc} below) that this module
is irreducible if the differences $\la_i-\mu_j$ are not integers.
Furthermore, for any $c\in\CC$ the simultaneous shifts 
$\la_i\mapsto \la_i+c$ and $\mu_j\mapsto \mu_j+c$ for all $i$ and $j$
do not affect the irreducibility of $L(\la)\ot L(\mu)$; see
Proposition~\ref{prop:shift}.
Thus, we may assume without loss of generality that all the entries
of $\la$ and $\mu$ are integers.

We shall be using the following definition.
Two disjoint finite subsets $A$ and $B$ of $\ZZ$ are {\it crossing\/} 
if there exist elements $a_1,a_2\in A$ and $b_1,b_2\in B$ such that
either $a_1<b_1<a_2<b_2$, or $b_1<a_1<b_2<a_2$. Otherwise, $A$ and $B$
are called {\it non-crossing\/}.

Given a highest weight $\la$ with integer entries introduce the
following subset of $\ZZ$:
\beql{alambda}
\A_{\la}=\{\la_1,\la_2-1,\dots,\la_n-n+1\}.
\non
\end{equation}

\bth\label{thm:main1} 
The module $L(\la)\ot L(\mu)$ is irreducible if and only if
the sets $\A_{\la}\setminus \A_{\mu}$ and $\A_{\mu}\setminus \A_{\la}$
are non-crossing.
\eth

Using the argument of Kitanine, Maillet and Terras~\cite{kmt:ff, mt:oq},
Nazarov and Tarasov \cite[Theorem~4.9]{nt:it} demonstrated that the irreducibility
criterion for the multiple tensor product \eqref{tenev} can be obtained
from the particular case of $k=2$ tensor factors.
Namely, the following ``binary property" holds. Here we
let $\la^{(1)},\dots,\la^{(k)}$ be $n$-tuples 
of complex numbers such that
$\la^{(p)}_i-\la^{(p)}_{i+1}$ is a non-negative integer for each $i$ and $p$.

\bth\label{thm:main2} The module
\beql{tenev0}
L(\la^{(1)})\ot L(\la^{(2)})\ot\cdots\ot L(\la^{(k)})
\end{equation}
is irreducible
if and only if all the modules $L(\la^{(p)})\ot L(\la^{(q)})$ 
with $p< q$ are irreducible.
\eth

Note that the ``only if" part of this theorem is well known.  It is implied by
Proposition~\ref{prop:perm} (see below); cf. \cite{cp:yr}, \cite{m:fd}.
If the module \eqref{tenev0} is irreducible, its highest weight
is easy to find. Therefore, together with Theorem~\ref{thm:main1}
the binary property of Theorem~\ref{thm:main2} allows one to determine
whether a given irreducible $\Y(\gl_n)$-module can be realized in a tensor product
\eqref{tenev0}.

For the proof of Theorem~\ref{thm:main1} we use 
the Gelfand--Tsetlin
bases of the $\gl_n$-modules $L(\la)$ and $L(\mu)$.
The key role is played by the formulas for the action
of the Drinfeld generators of $\Y(\gl_n)$ in these bases
as well as by the quantum minor formulas for the Yangian lowering operators
\cite{m:gt}; cf. Nazarov and Tarasov~\cite{nt:yg, nt:ry}.

In the case of the Yangian $\Y(\gl_2)$ the criterion coincides with the one
obtained by Chari and Pressley~\cite{cp:yr}	
and it is also implicitly contained
in Tarasov's paper~\cite{t:im}; see also \cite{m:fd}.
Nazarov and Tarasov~\cite{nt:tp} 
found a criterion of irreducibility of \eqref{tenev}
in the case where
each highest weight $\lambda^{(p)}$
has the form $(\al,\dots,\al,\be,\dots,\be)$ with $\al-\be\in\ZZ_+$.
This generalized earlier results by Akasaka and Kashiwara~\cite{ak:fd}
and Zelevinsky~\cite{z:ir}.

Leclerc, Nazarov and Thibon~\cite{lnt:ir} have found an irreducibility criterion for
the induction products of evaluation modules over the affine Hecke
algebras of type ${\rm A}$ with the use of the canonical bases;
see also Leclerc and Thibon~\cite{lt:ri},
Leclerc and Zelevinsky~\cite{lz:qf}. 
The application of the Drinfeld functor 
\cite{d:da}  (see also \cite{a:df}) leads to
an irreducibility criterion for the Yangian modules \eqref{tenev}
(equivalent to Theorems~\ref{thm:main1} and \ref{thm:main2}), when the
highest weights $\la^{(p)}$ satisfy some extra conditions. Namely,
assuming that the $\la^{(p)}$ are partitions (we may do this without loss
of generality), one should require that the sum of their lengths does not
exceed $n$.
 
\medskip
This work was completed during the author's visit
to the Erwin Schr\"odinger Institute, Vienna
for the Representation Theory Program, 2000. 
I would like to thank ESI and
the organizers of the program, A.~A.~Kirillov and 
V.~G.~Kac for the invitation.
I~am very much obliged to M.~L.~Nazarov for valuable discussions
during the visit. 
He kindly informed me about the results of \cite{lnt:ir, nt:it} prior to their
publication, and the present form of Theorem~\ref{thm:main1} was inspired
by the irreducibility criterion in \cite{lnt:ir}. Originally,
the author obtained this result in the form
of Theorems~\ref{thm:sufc} and \ref{thm:necc}: see below.
The financial support from the Australian Research Council is
acknowledged.

\section{Preliminaries}\label{sec:pre}
\setcounter{equation}{0}

We refer the reader to the expository papers \cite{m:fd,mno:yc}
where the results on the structure theory and representations
of the Yangians are collected.

The {\it Yangian\/} $\Y(n)=\Y(\gl_n)$ \cite{d:ha, d:nr} is the
complex associative algebra with the
generators $t_{ij}^{(1)},t_{ij}^{(2)},\dots$ where $1\leq i,j\leq n$,
and the defining relations
\beql{defrel}
[t_{ij}(u),t_{kl}(v)]=\frac{1}{u-v}\Big(t_{kj}(u)t_{il}(v)-t_{kj}(v)t_{il}(u)\Big),
\end{equation}
where
\beql{series}
t_{ij}(u) = \delta_{ij} + t^{(1)}_{ij} u^{-1} + t^{(2)}_{ij}u^{-2} +
\cdots \in \Y(n)[[u^{-1}]]
\non
\end{equation}
and $u$ is a formal (commutative) variable. 
The Yangian $\Y(n)$ is a Hopf algebra with the coproduct 
$
\Delta : \Y(n)\to
\Y(n)\ot\Y(n) 
$ 
defined by
\beql{copr}
\Delta (t_{ij}(u))=\sum_{a=1}^n t_{ia}(u)\ot
t_{aj}(u).
\end{equation}

Given sequences $a_1,\ldots, a_r$ and $b_1,\ldots, b_r$ of elements
of $\{1,\dots,n\}$ the corresponding {\it quantum minor\/} 
of the matrix $\big[t_{ij}(u)\big]$
is defined
by the following equivalent formulas:
\begin{align}\label{qminor}
{t\ts}^{a_1\cdots\ts a_r}_{b_1\cdots\ts b_r}(u)&=
\sum_{\sigma\in \Sym_r} \sgn \sigma \cdot t_{a_{\sigma(1)}b_1}(u)\cdots
t_{a_{\sigma(r)}b_r}(u-r+1)\\
\label{qminor2}
{}&=
\sum_{\sigma\in \Sym_r} \sgn \sigma\cdot t_{a_1b_{\sigma(1)}}(u-r+1)\cdots
t_{a_rb_{\sigma(r)}}(u).
\end{align}
The series ${t\ts}^{a_1\cdots\ts a_r}_{b_1\cdots\ts b_r}(u)$ is skew symmetric
under permutations of the indices $a_i$, or $b_i$.

The Poincar\'e--Birkhoff--Witt theorem for the Yangian $\Y(n)$ (see e.g. 
\cite[Corollary~1.23]{mno:yc}) implies that given a subset of indices
$\{a_1,\dots,a_r\}\subseteq \{1,\dots,n\}$ the coefficients of the series
$t_{a_ia_j}(u)$ with $i,j=1,\dots,r$ generate a subalgebra of $\Y(n)$
isomorphic to $\Y(r)$.

The mapping
\beql{audual}
t_{ij}(u)\mapsto 
(-1)^{i+j}\ts{t\ts}^{1\cdots\ts\wh j\ts\cdots n}_{1\cdots\ts\wh i\ts\cdots\ts n}(-u)
\end{equation}
defines an algebra automorphism of $\Y(n)$; the hats indicate the indices
to be omitted.

The following proposition is proved in
\cite[Proposition~1.1]{m:yt} by using
the $R$-matrix form of the defining relations \eqref{defrel}.

\bpr\label{prop:qmr}
We have the relations
\beq
\begin{aligned}
{}&[{t\ts}^{a_1\cdots\ts a_k}_{b_1\cdots\ts b_k}(u),
{t\ts}^{c_1\cdots\ts c_l}_{d_1\cdots\ts d_l}(v)]=
\sum_{p=1}^{\min\{k,l\}}\frac{(-1)^{p-1}\ts p!}
{(u-v-k+1)\cdots (u-v-k+p)}
\\
{}&\sum_{\overset{\scl i_1<\cdots<i_p}
{\scl j_1<\cdots<j_p}}\left(
{t\ts}^{a_1\cdots\ts c_{j_1}\cdots\ts 
c_{j_p}\cdots\ts a_k}_{b_1\ \cdots\ \  b_k}(u)
{t\ts}^{c_1\cdots\ts a_{i_1}\cdots\ts 
a_{i_p}\cdots\ts c_l}_{d_1\ \cdots\ \  d_l}(v)
-{t\ts}^{c_1\ \cdots\ \  c_l}_{d_1\cdots\ts 
b_{i_1}\cdots\ts b_{i_p}\cdots\ts d_l}(v)
{t\ts}^{a_1\ \cdots\ \  a_k}_{b_1\cdots\ts 
d_{j_1}\cdots\ts d_{j_p}\cdots\ts b_k}(u)
\right).
\end{aligned}
\non
\end{equation}
Here the $p$-tuples
of upper indices $(a_{i_1},\dots, a_{i_p})$ and $(c_{j_1},\dots, c_{j_p})$
are respectively
interchanged in the first summand on the right hand
side while the $p$-tuples of lower indices
$(b_{i_1},\dots, b_{i_p})$ and
$(d_{j_1},\dots, d_{j_p})$
are interchanged in the second
summand.  \endproof
\epr

We note the following particular case of these relations:
\beql{qmrel}
[{t}^{}_{ab}(u),
{t\ts}^{c_1\cdots\ts c_l}_{d_1\cdots\ts d_l}(v)]=
\frac{1}
{u-v}
\left(\sum_{i=1}^l
{t}^{}_{c_ib}(u)
{t\ts}^{c_1\cdots\ts a\ts\cdots\ts c_l}_{d_1\ \cdots\ \  d_l}(v)
-\sum_{i=1}^l{t\ts}^{c_1\ \cdots\ \  c_l}_{d_1\cdots\ts 
b\ts\cdots\ts d_l}(v)
{t}^{}_{ad_i}(u)\right).
\end{equation}
This implies the well-known 
property of the quantum minors: for any indices $i,j$ we have
\beql{center}
[{t}^{}_{c_id_j}(u),
{t\ts}^{c_1\cdots\ts c_l}_{d_1\cdots\ts d_l}(v)]=0.
\end{equation}

We shall frequently use the following result proved in \cite{nt:ry}.

\bpr\label{prop:delqm}
The images of the quantum minors under the coproduct
are given by
\beql{delqm}
\Delta ({t\ts}^{a_1\cdots\ts a_r}_{b_1\cdots\ts b_r}(u))
=\sum_{c_1<\cdots < c_r}{t\ts}^{a_1\cdots\ts a_r}_{c_1\cdots\ts c_r}(u)
\ot {t\ts}^{c_1\cdots\ts c_r}_{b_1\cdots\ts b_r}(u),
\non
\end{equation}
summed over all subsets of indices $\{c_1,\dots,c_r\}$ from $\{1,\dots,n\}$.
\endproof
\epr

For $m\geq 1$ introduce the series $a_m(u)$, $b_m(u)$ and $c_m(u)$ by
\beql{drgen}
a_m(u)={t\ts}^{1\cdots\ts m}_{1\cdots\ts m}(u),\qquad
b_m(u)={t\ts}^{1\cdots\ts m}_{1\cdots\ts m-1,m+1}(u), \qquad
c_m(u)={t\ts}^{1\cdots\ts m-1,m+1}_{1\cdots\ts m}(u).
\end{equation}
The coefficients of these series generate the algebra $\Y(n)$ \cite{d:nr},
they are called the {\it Drinfeld generators\/}.

By a theorem of Drinfeld \cite{d:nr} every finite-dimensional irreducible
representation of the Yangian $\Y(n)$ is 
a highest weight representation. That is,
it contains a unique, up to a scalar factor, nonzero vector $\zeta$
(the {\it highest vector\/})
which is annihilated by all upper triangular elements $t_{ij}(u)$, $i<j$,
and $\zeta$ is an eigenvector for the diagonal generators $t_{ii}(u)$,
\beq
t_{ii}(u)\ts \zeta=\lambda_i(u)\ts\zeta,\qquad i=1,\dots,n.
\non
\end{equation}
Here the $\la_i(u)$ are formal series in $u^{-1}$ with complex coefficients.
We call the collection $(\la_1(u),\dots,\la_n(u))$ the {\it highest weight\/}
of the representation. 
Equivalently, $\zeta$ is annihilated by $b_1(u),\dots,b_{n-1}(u)$
and it is an eigenvector for each of the operators $a_1(u),\dots,a_n(u)$
\cite{d:nr}, so that
\beq
a_{m}(u)\ts \zeta=\lambda_1(u)\lambda_2(u-1)\cdots
\lambda_m(u-m+1)\ts\zeta,\qquad m=1,\dots,n.
\non
\end{equation}

If $L$ is any $\Y(n)$-module, then a nonzero
element $\zeta\in L$ is called a {\it singular vector\/}
if $\zeta$ is annihilated by all upper triangular generators $t_{ij}(u)$, $i<j$,
and $\zeta$ is an eigenvector for the diagonal elements $t_{ii}(u)$.
Such a vector $\zeta$
generates a highest weight submodule in $L$.
The following proposition if proved by a standard argument;
see e.g. \cite{m:fd}.

\bpr\label{prop:sing}
If $L$ is an irreducible 
highest weight $\Y(n)$-module and $\zeta\in L$
is annihilated by all operators $t_{ij}(u)$ with $i<j$ then
$\zeta$ is proportional to the highest vector of $L$.
\endproof
\epr

Let the $E_{ij}$, $i,j=1,\dots,n$ denote the standard basis elements 
of the Lie algebra $\gl_n$. For any $a\in\CC$ the mapping 
\beql{epi}
\varphi_a:t_{ij}(u)\mapsto \delta_{ij}+\frac{E_{ij}}{u-a}
\end{equation}
defines an algebra epimorphism from $\Y(n)$ to the universal enveloping algebra
$\U(\gl_n)$ so that any $\gl_n$-module can be extended to a $\Y(n)$
module via \eqref{epi}. In particular, let $\la$
be an $n$-tuple
of complex numbers $\la=(\la_1,\dots,\la_n)$ such that $\la_i-\la_{i+1}\in\ZZ_+$
for all $i$ (we call such $n$-tuples $\gl_n$-{\it highest weights\/}).
Consider the irreducible finite-dimensional
$\gl_n$-module $L(\la)$
with the highest weight $\la$
with respect to the upper triangular Borel subalgebra.
The corresponding $\Y(n)$-module is denoted by
$L_a(\la)$, and we call	it
the {\it evaluation module\/}. We keep the notation
$L(\la)$ for the module $L_a(\la)$ with $a=0$.
The coproduct $\Delta$ defined by \eqref{copr} 
allows one to consider the tensor products
\eqref{tenev} as $\Y(n)$-modules.

Let us denote by $I$ the $n$-tuple $(1,1,\dots,1)$.

\bpr\label{prop:redu} The $\Y(n)$-module \eqref{tenev} is irreducible
if and only if the module
\beql{tensh}
L(\la^{(1)}-a_1\ts I)\ot L(\la^{(2)}-a_2\ts I)
\ot\cdots\ot L(\la^{(k)}-a_k\ts I)
\end{equation}
is irreducible.
\epr

\Proof Suppose that \eqref{tenev} is irreducible. Let $\xi^{(p)}$ denote
the highest vector of the $\gl_n$-module $L(\la^{(p)})$. 
We derive from \eqref{copr} and \eqref{epi} that
$\zeta=\xi^{(1)}\ot\cdots\ot \xi^{(k)}$ is the highest vector
of the $\Y(n)$-module \eqref{tenev} with the highest weight
$(\la_1(u),\dots,\la_n(u))$ where
\beql{tenhw}
\la_i(u)=\Big(1+\frac{\la_i^{(1)}}{u-a_1}\Big)\cdots 	
\Big(1+\frac{\la_i^{(k)}}{u-a_k}\Big),\qquad
i=1,\dots,n.
\end{equation}
Consider the automorphism of $\Y(n)$
\beql{scal}
t_{ij}(u)\mapsto f(u)\ts t_{ij}(u),
\end{equation}
where $f(u)$ is the formal series in $u^{-1}$ given by
\beql{fu}
f(u)=(1-a_1u^{-1})\cdots (1-a_ku^{-1}).
\non
\end{equation}
The composition of the module \eqref{tenev} with this automorphism
is an irreducible $\Y(n)$-module $\wt{L}$ with the highest weight
$(\wt{\la}_1(u),\dots,\wt{\la}_n(u))$ where
\beql{tenhwwt}
\wt{\la}_i(u)=\Big(1+\frac{\la_i^{(1)}-a_1}{u}\Big)\cdots 	
\Big(1+\frac{\la_i^{(k)}-a_k}{u}\Big),\qquad
i=1,\dots,n.
\end{equation}
On the other hand, the tensor product of the highest vectors
of the $\gl_n$-modules $L(\la^{(p)}-a_p\ts I)$
is a singular vector of the $\Y(n)$-module \eqref{tensh}
with the weight given by \eqref{tenhwwt}. Therefore, $\wt{L}$
is isomorphic to a subquotient of \eqref{tensh}. However,
these two modules have the same dimension and hence, they are isomorphic.
In particular, the module \eqref{tensh} is irreducible.
The proof is completed by reversing the argument.
\endproof

\bpr\label{prop:shift} Given $c\in\CC$, the simultaneous shifts
\beql{shift}
\la^{(p)}_i\mapsto \la^{(p)}_i+c,\qquad p=1,\dots,k,\quad i=1,\dots n
\non
\end{equation}
of the parameters of the module \eqref{tenev0} do not affect its
irreducibility.
\epr

\Proof	It the module \eqref{tenev0} is irreducible then so is 
the module $L_{}^{c}$ which is the composition of \eqref{tenev0}
with the automorphism of $\Y(n)$ given by
\beql{shiauto}
t_{ij}(u)\mapsto t_{ij}(u+c),\qquad i,j=1,\dots,n.
\non
\end{equation}
The highest weight of $L_{}^{c}$ is $(\la^c_1(u),\dots,\la^c_n(u))$ with
\beql{tenhwc}
\la^c_i(u)=\Big(1+\frac{\la_i^{(1)}}{u+c}\Big)\cdots 	
\Big(1+\frac{\la_i^{(k)}}{u+c}\Big),\qquad
i=1,\dots,n.
\non
\end{equation}
The proof is completed by repeating the argument of the proof
of Proposition~\ref{prop:redu} with the use of 
the automorphism \eqref{scal}
of $\Y(n)$ where
$f(u)=(1+cu^{-1})^k$.			
\endproof

\bpr\label{prop:perm}
Suppose that the $\Y(n)$-module \eqref{tenev0} is irreducible.
Then any permutation of the tensor factors in \eqref{tenev0}
gives an isomorphic representation of $\Y(n)$.
\epr

\Proof
Denote the tensor product \eqref{tenev0} by $L$. Note that
$L$ is a representation with the highest weight
$(\la_1(u),\dots,\la_n(u))$
given by 
\eqref{tenhw} with $a_1=\cdots=a_k=0$.
Consider a representation $L'$ obtained
by a certain permutation of the tensor factors in \eqref{tenev0}.
The tensor product $\zeta'$ of the highest vectors
of the representations $L(\lambda^{(i)})$ is a singular
vector in $L'$ whose weight is given by the same formulas \eqref{tenhw}.
This implies that $\zeta'$ generates a highest weight
submodule in $L'$ such that	its irreducible quotient is isomorphic
to $L$. However, $L$ and $L'$ have the same dimension which implies that
$L$ and $L'$ are isomorphic. 
\endproof

We shall use a version given in \cite{m:gt} (cf. \cite{nt:yg})
of the construction of a basis of the $\gl_n$-module $L(\la)$
which is originally due to Gelfand and Tsetlin \cite{gt:fd}.
We equip $L(\la)$ with a $\Y(n)$-module structure by using
the epimorphism
\beql{epi0}
\Y(n)\to\U(\gl_n),\qquad 
t_{ij}(u)\mapsto \delta_{ij}+E_{ij}\ts u^{-1},
\end{equation}
see \eqref{epi}.
A {\it pattern\/} $\La$ (associated with $\la$) is a sequence
of rows $\La_n,\La_{n-1},\dots,\La_1$, where
$\La_r=(\lambda_{r1},\dots,\lambda_{rr})$ is the $r$-th row from the bottom, 
the top row $\La_n$
coincides with $\la$, and the following
{\it betweenness conditions\/} are satisfied: for $r=2,\dots,n$
\beql{betw}
\la_{ri}-\la_{r-1,i}\in\ZZ_+,\quad 
\la_{r-1,i}-\la_{r,i+1}\in\ZZ_+,\qquad{\rm for}\quad i=1,\dots,r-1.
\end{equation}
For any pattern $\La$ introduce the vector $\xiL\in L(\la)$ by
\beql{xildef}
\xiL=\prod_{r=2,\dots,n}^{\rightarrow}\prod_{i=1}^{r-1}
\tau^{}_{ri}(-\la_{r-1,i}-1)\cdots \tau^{}_{ri}(-\la_{ri}+1)
\tau^{}_{ri}(-\la_{ri})\ts \xi,
\end{equation}
where $\xi$ is the highest vector of $L(\la)$ and
\beql{taugln}
\tau^{}_{ri}(u)=u(u-1)\cdots (u-r+i+1)\ts 
{t\ts}^{i+1\ts\cdots\ts r}_{i\ts\cdots\ts r-1}(u)
\non
\end{equation}
is the {\it lowering operator\/}; see also Section~\ref{sec:necc}.
The vectors $\xiL$, where $\La$ runs over all patterns associated with $\la$,
form a basis of $L(\la)$. The $\tau^{}_{ri}(u)$ essentially
coincide with the standard lowering operators arising from
the transvector algebras; cf. \cite{m:yt}.
We find from \eqref{epi0} that the operators 
\beq
\bal
B_m(u)&=u(u-1)\cdots (u-m+1)\ts b_m(u),\\
A_m(u)&=u(u-1)\cdots (u-m+1)\ts a_m(u)
\eal
\non
\end{equation}
in $L(\la)$ are polynomials in $u$; see \eqref{drgen}.
Their action
in the basis $\{\xiL\}$ of $L(\la)$ 
is given by the following formulas; see \cite{m:gt}.
They can also be deduced from Lemmas~\ref{lem:tii+1}--\ref{lem:yn-1};
see Section~\ref{sec:necc}.
We use the notation	 $l_{ri}=\la_{ri}-i+1$.

\bpr\label{prop:agt} 
We have
\begin{align}\label{agt}
A_m(u)\ts\xiL&=(u+l_{m1})\cdots (u+l_{mm})\ts\xiL,\\
\label{bgt}
B_m(-l_{mj})\ts\xiL&=-\prod_{i=1}^{m+1} (l_{m+1,i}-l_{mj})
\ts \xi^{}_{\La+\delta_{mj}}\quad\text{for}\quad j=1,\dots,m,
\non
\end{align}					   
where $\La+\delta_{mj}$ is obtained from
$\La$ by replacing the entry $\lambda_{mj}$ with $\lambda_{mj}+1$,
and $\xi^{}_{\La+\delta_{mj}}$ is supposed to be equal to zero
if $\La+\delta_{mj}$ is not a pattern. \endproof
\epr

Applying the Lagrange interpolation formula we can find 
the action
of $B_m(u)$ for any $u$. Note that the polynomial $B_m(u)$	has degree $m-1$
with the leading coefficient $E_{m,m+1}$. This therefore implies 
the Gelfand--Tsetlin formulas \cite{gt:fd} for the action of 
the elements $E_{m,m+1}$:
\beql{Egt}
E_{m,m+1}\ts\xiL=-\sum_{j=1}^m\frac{(l_{m+1,1}-l_{mj})\cdots (l_{m+1,m+1}-l_{mj})}
{(l_{m1}-l_{mj})\cdots\wedge_j\cdots (l_{mm}-l_{mj})}
\ts \xi^{}_{\La+\delta_{mj}},
\end{equation}
where $\wedge_j$ indicates that the $j$-th factor is skipped.

We conclude this section with an equivalent form of the conditions of
Theorem~\ref{thm:main1}. 

Given complex $\gl_n$-highest weights $\la=(\la_1,\dots,\la_n)$ and 
$\mu=(\mu_1,\dots,\mu_n)$ we shall use the notation
\beql{lm}
l_i=\la_i-i+1,\quad m_i=\mu_i-i+1,\qquad i=1,\dots,n.
\non
\end{equation}
For a pair of indices $i<j$ we shall denote
\beql{lanran} 
\bal
\lan l_j,l_i\ran&=\{l_j,l_j+1,\dots,l_i\}\setminus\{l_j,l_{j-1},\dots,l_i\},	
\\
\lan m_j,m_i\ran&=\{m_j,m_j+1,\dots,m_i\}\setminus\{m_j,m_{j-1},\dots,m_i\}.
\eal
\end{equation}
In particular, if $\la_i=\la_{i+1}=\cdots=\la_j$ then 
$\lan l_j,l_i\ran =\emptyset$.

We shall assume now 
that $\la$ and $\mu$ are $\gl_n$-highest weights
with integer entries.

\bpr\label{prop:equiv}
The sets $\A_{\la}\setminus \A_{\mu}$ and $\A_{\mu}\setminus \A_{\la}$
are non-crossing if and only if for all pairs of indices $1\leq i< j\leq n$
we have
\beql{ijcondeq}
m_j,m_i\not\in\lan l_j,l_i\ran\qquad\text{or}
\qquad l_j,l_i\not\in\lan m_j,m_i\ran.
\end{equation}
\epr

\Proof  Let us write $\Cond(\A_{\la},\A_{\mu})$ for the condition
that $\A_{\la}\setminus \A_{\mu}$ and $\A_{\mu}\setminus \A_{\la}$
are non-crossing.
We use induction on $n$. In the case $n=2$ the statement is obviously true.
Let $n\geq 3$. Suppose first that \eqref{ijcondeq} holds.
Set
\beql{alapm}
\A_{\la}^-=\{l_1,\dots,l_{n-1}\}\qquad\text{and}\qquad
\A_{\la}^+=\{l_2,\dots,l_{n}\}
\non
\end{equation}
and similarly define $\A_{\mu}^-$ and $\A_{\mu}^+$.
By the induction hypothesis, both conditions
$\Cond(\A_{\la}^-,\A_{\mu}^-)$ and $\Cond(\A_{\la}^+,\A_{\mu}^+)$
are satisfied.
If $m_1=l_1$ then $\Cond(\A_{\la},\A_{\mu})$ obviously holds. 
We may assume without loss of generality that $m_1>l_1$.
Let
\beql{almam}
\A_{\la}^-\setminus \A_{\mu}^-=\{l_{i_1},\dots,l_{i_k}\},
\qquad
\A_{\mu}^-\setminus \A_{\la}^-=\{m_{j_1},\dots,m_{j_k}\},
\non
\end{equation}
where $1\leq i_1<\cdots<i_k\leq n$ and $1=j_1<\cdots<j_k\leq n$.
We must have for some $a\in\{1,\dots,k\}$ that
\beql{lbetw}
m_{j_{a+1}} < l_{i_{k}}< \cdots <	l_{i_1}< m_{j_a},
\non
\end{equation}
where the leftmost inequality is ignored when $a=k$.
If $2\leq a\leq k-1$ then together with $\Cond(\A_{\la}^+,\A_{\mu}^+)$
this clearly ensures $\Cond(\A_{\la},\A_{\mu})$.
Similarly, this is also true when $a=1$ and $i_1\geq 2$.
So, if $a=1$ then the only case where both 
$\Cond(\A_{\la}^-,\A_{\mu}^-)$ and $\Cond(\A_{\la}^+,\A_{\mu}^+)$ hold
but $\Cond(\A_{\la},\A_{\mu})$ does not, is the one with the following
inequalities between the elements of $\A_{\la}\setminus \A_{\mu}$
and $\A_{\mu}\setminus \A_{\la}$:
\beql{pmh}
l_n<m_n<m_{j_k}<\cdots<m_{j_2}<l_{i_k}<\cdots<l_{i_2}<l_1<m_1.
\non
\end{equation}
However, in this case $m_n\in\lan l_n,l_1\ran$ and $l_1\in \lan m_n,m_1\ran$,
so that \eqref{ijcondeq} is violated
for $i=1$ and $j=n$.
An analogous argument shows that if $a=k$ then the only case where both 
$\Cond(\A_{\la}^-,\A_{\mu}^-)$ and $\Cond(\A_{\la}^+,\A_{\mu}^+)$ hold
but $\Cond(\A_{\la},\A_{\mu})$ does not, is 
\beql{pmh2}
l_n<l_{i_k}<\cdots<l_{i_2}<m_n<l_1<m_{j_k}<\cdots<m_{j_2}<m_1.
\non
\end{equation}
But then $m_n\in\lan l_n,l_1\ran$ and $l_1\in \lan m_n,m_1\ran$
which contradicts \eqref{ijcondeq} again.

Conversely, suppose that $\Cond(\A_{\la},\A_{\mu})$ holds. This condition clearly
implies both $\Cond(\A_{\la}^-,\A_{\mu}^-)$ and $\Cond(\A_{\la}^+,\A_{\mu}^+)$
and so, by the induction hypothesis, \eqref{ijcondeq} holds for all
pairs $i<j$ with the possible
exception for $(i,j)=(1,n)$. If the latter condition fails then we have
\beql{poss}
\left\{
\bal
m_n&\in \lan l_n,l_1\ran\\
l_1&\in \lan m_n,m_1\ran
\eal
\right.
\qquad\text{or}\qquad
\left\{
\bal
l_n&\in \lan m_n,m_1\ran\\
m_1&\in \lan l_n,l_1\ran.
\eal
\right.
\end{equation}
However, in each of the two cases this contradicts
$\Cond(\A_{\la},\A_{\mu})$. \endproof

\section{Sufficient conditions}\label{sec:sufc}
\setcounter{equation}{0}

Our aim in this section is to prove the following.

\bth\label{thm:sufc}
Let $\la$ and $\mu$ be complex $\gl_n$-highest weights. Suppose that
for each pair of indices $1\leq i<j\leq n$ we have
\beql{ijcond}
m_j,m_i\not\in\lan l_j,l_i\ran\qquad\text{or}
\qquad l_j,l_i\not\in\lan m_j,m_i\ran.
\end{equation}
Then the $\Y(n)$-module $L(\la)\ot L(\mu)$ is irreducible.
\eth

We give the proof of Theorem~\ref{thm:sufc} as a sequence of lemmas.
Let $\xi$ and $\xi'$ denote the highest vectors of the $\gl_n$-modules
$L(\la)$ and $L(\mu)$, respectively.
Let $N$ be a nonzero $\Y(n)$-submodule of $L(\la)\ot L(\mu)$.
A standard argument (see e.g. \cite{m:fd}) shows that $N$ must
contain a singular vector $\zeta$. 
The key part of the proof of the theorem is to show by induction 
on $n$ that
\beql{sing}
\zeta={\rm const}\cdot\xi\ot\xi'.
\end{equation}
Then considering dual modules we also show that
the vector $\xi\ot\xi'$ is cyclic.

By Proposition~\ref{prop:perm}, exchanging $\la$ and $\mu$ if necessary, 
we may assume that for the pair $(i,j)$ with
$i=1$ and $j=n$ the condition
\beql{1ncond}
m_1,m_n\not\in\lan l_n,l_1\ran
\end{equation}
is satisfied.
Consider
the Gelfand--Tsetlin basis $\{\xiL\}$ of the
$\gl_n$-module $L(\la)$; see Section~\ref{sec:pre}. The singular vector
$\zeta$ is uniquely written in the form
\beql{zeta}
\zeta=\sum_{\La}\xiL\ot \etaL,
\end{equation}
summed over all patterns $\La$ associated with $\la$, and
$\etaL\in L(\mu)$.

For the diagonal Cartan subalgebra $\mathfrak h$ of $\gl_n$,
we denote by $\ve_i$  the basis vector	of $\mathfrak h^*$
dual to the element $E_{ii}$ so that the $n$-tuple $\la$ can be identified
with the element $\la_1\ve_1+\cdots+\la_n\ve_n\in \mathfrak h^*$.
We shall be using a standard partial ordering on the weights of $L(\la)$.
Given two weights $v,w\in\mathfrak h^*$,
we shall write $v\preceq w$ if $w-v$ is a
$\ZZ_+$-linear combination
of the simple roots
$\ve_a-\ve_{a+1}$. Equivalently, $v\preceq w$ if and only if
\beql{eqp}
w-v=\sum_{a=1}^{n}p_a\ve_a, 
\end{equation}
with the conditions
\beql{eqpco}
p_1,\ p_1+p_2,\ \dots,\ p_1+\cdots+p_{n-1}\in\ZZ_+,\quad 
p_1+\cdots+p_{n}=0.
\non
\end{equation}

The embedding 
\beql{emb}
\U(\gl_n)\hookrightarrow\Y(n),\qquad E_{ij}\mapsto t_{ij}^{(1)},
\end{equation}
defines the natural
$\U(\gl_n)$-module  structure
on $L(\la)\ot L(\mu)$. We shall usually identify the operators 
$E_{ij}$ and $t_{ij}^{(1)}$.
The vector $\zeta$ is clearly a $\gl_n$-singular vector. In particular,
it is a weight vector. Since the basis $\{\xiL\}$ consists of 
weight vectors, each element $\etaL\in L(\mu)$ in \eqref{zeta} is also
a $\gl_n$-weight vector. Moreover, all elements $\xiL\ot\etaL$ in \eqref{zeta} have
the same $\gl_n$-weight.

We shall denote the weight of the vector $\xiL$, or, the weight
of the pattern $\La$, by $w(\La)$.
It is well known~\cite{gt:fd}, and can be deduced e.g. from \eqref{agt} that
\beql{weigt}
w(\La)=w_1\ve_1+\cdots+w_n\ve_n,\qquad w_k=
\sum_{i=1}^k\la_{ki}-\sum_{i=1}^{k-1}\la_{k-1,i}.
\end{equation}

We shall say that a pattern $\La$ {\it occurs\/} in the expansion \eqref{zeta}
if $\eta^{}_{\La}\ne 0$. Consider the set of patterns occurring in \eqref{zeta}
and suppose that $\La^0$ is a minimal element of this set with respect to the
partial ordering on the weights $w(\La)$. In other words, if $\La$ occurs in \eqref{zeta}
and $w(\La)\preceq w(\La^0)$
then $w(\La)=w(\La^0)$.

\ble\label{lem:lao}
The vector
$\eta^{}_{\La^0}$ coincides with $\xi'$, up to a constant factor. 
\ele

\Proof
We have $b_m(u)\ts\zeta=0$ for $m=1,\dots,n-1$. Therefore, by
Proposition~\ref{prop:delqm},
\beql{bzet}
\sum_{c_1<\cdots<c_m}\sum_{\La}
{t\ts}^{1\cdots\ts m}_{c_1\cdots\ts c_m}(u)\ts\xiL\ot
{t\ts}^{c_1\cdots\ts c_m}_{1\cdots\ts m-1,m+1}(u)\ts\etaL=0.
\end{equation}
Write the elements ${t\ts}^{1\cdots\ts m}_{c_1\cdots\ts c_m}(u)\ts\xiL$
as linear combinations of the basis vectors $\xiL$ and take the coefficient
at $\xiLo$ in the relation \eqref{bzet}.
The weight of the vector ${t\ts}^{1\cdots\ts m}_{c_1\cdots\ts c_m}(u)\ts\xiL$
(or, to be more precise, the weight of each of the coefficients of these series)
equals 
\beq
w'=w(\La)+\ve_1+\cdots+\ve_m-\ve_{c_1}-\cdots-\ve_{c_m}.
\non
\end{equation}
Therefore $w'\succeq w(\La)$, and $w'=w(\La)$ if and only if $c_i=i$ for each $i$.
Since $\La^0$ is a pattern of a minimal weight, the vector $\xiLo$
can only occur in the expansion of 
\beql{expa}
{t\ts}^{1\cdots\ts m}_{1\cdots\ts m}(u)\ts\xiLo=a_m(u)\ts\xiLo.
\end{equation}
By \eqref{agt} this implies that
\beql{etaL}
{t\ts}^{1\cdots\ts m}_{1\cdots\ts m-1,m+1}(u)\ts\etaLo=b_m(u)\ts\etaLo=0
\non
\end{equation}
for each $m$. Thus, $\etaLo$ is a singular vector of $L(\mu)$
and so,
$\etaLo={\rm const}\cdot\xi'$.	
\endproof

\ble\label{lem:onep} 
The pattern $\La^0$ is determined uniquely.
Moreover, if a pattern $\La$ occurs in \eqref{zeta}
then $w(\La)\succeq w(\La^0)$.
\ele

\Proof
We use the fact that $\zeta$ is an eigenvector
for the generators $a_1(u),\dots,a_n(u)$. By
Proposition~\ref{prop:delqm},
\beql{amzeta}
a_m(u)\ts\zeta=\sum_{c_1<\cdots<c_m}\sum_{\La}
{t\ts}^{1\cdots\ts m}_{c_1\cdots\ts c_m}(u)\ts\xiL\ot
{t\ts}^{c_1\cdots\ts c_m}_{1\cdots\ts m}(u)\ts\etaL.
\end{equation}
This vector equals $\al_m(u)\ts	\zeta$ for a formal series $\al_m(u)$.
As we have seen in the proof of Lemma~\ref{lem:lao},
the vector $\xiLo$ can only occur 
in \eqref{amzeta} in the expansion of \eqref{expa}.
By \eqref{agt} and Lemma~\ref{lem:lao},
comparing the coefficients at $\xiLo$, we find that
\beql{comp}
\al_m(u)=\frac{(u+l^{\ts 0}_{m1})\cdots (u+l^{\ts 0}_{mm})}
{u(u-1)\cdots (u-m+1)}\ts\be_m(u), 
\end{equation}
where $l^{\ts 0}_{mi}=\la^{\ts 0}_{mi}-i+1$ and the series $\be_m(u)$ is defined by
$a_m(u)\ts\xi'=\be_m(u)\ts\xi'$. The equation \eqref{comp} uniquely
determines the parameters $l^{\ts 0}_{mi}$, $1\leq i\leq m<n$, since
$l^{\ts 0}_{mi}-l^{\ts 0}_{m,i+1}$ is a positive integer for each $i$.
Thus, the pattern $\La^0$ is also determined uniquely.
The second claim is now obvious.
\endproof

\ble\label{lem:dif}
For each entry $\la^{}_{ka}$ of a pattern $\La$ 
occurring in \eqref{zeta} we have
\beql{lad}
\la^{}_{ka}-\la^0_{ka}\in\ZZ_+,\qquad\text{for}\quad 1\leq a\leq k\leq n-1.
\non
\end{equation}
\ele

\Proof
If $\La$ occurs in \eqref{zeta} 
then $w(\La)\succeq w(\La^0)$ by Lemma~\ref{lem:onep}.
We use induction on $w(\La)$. Fix $\La\ne\La^0$.
Then there exists another pattern $\La'$ occurring in \eqref{zeta}
such that $w(\La^0)\preceq w(\La')\prec w(\La)$, and 
for some $m$ and some indices $c_1<\cdots<c_m$ the expansion of
${t\ts}^{1\cdots\ts m}_{c_1\cdots\ts c_m}(u)\ts\xi^{}_{\La'}$ contains
$\xiL$ with a nonzero coefficient. Indeed, if this is not
the case, then considering the coefficient at $\xiL$
in \eqref{bzet}	we come to the conclusion	that
$b_m(u)\ts\etaL=0$ for all $m=1,\dots,n-1$, and so,
$\etaL$ is, up to a constant, the highest vector
of $L(\mu)$: see the proof of Lemma~\ref{lem:lao}. 
This implies that $\La$ and $\La^0$ must have the same weight.
Due to Lemma~\ref{lem:onep}, 
we have to conclude that $\La=\La^0$, contradiction.

By \eqref{qmrel} 
the operator ${t\ts}^{1\cdots\ts m}_{c_1\cdots\ts c_m}(u)$
can be represented as the commutator
\beql{tcommu}
{t\ts}^{1\cdots\ts m}_{c_1\cdots\ts c_m}(u)
=[\ldots[\ts[{t\ts}^{1\cdots\ts m}_{1\cdots\ts m}(u),
E_{mc_m}],\ts E_{{m-1},c_{m-1}}],
\dots,E_{pc_p}],
\non
\end{equation}
where $p$ is the minimum of the indices $i$ such that $c_i\ne i$.
Here, as before, we identify the elements $E_{ij}$ and $t_{ij}^{(1)}$
using the embedding \eqref{emb}.
The operator ${t\ts}^{1\cdots\ts m}_{1\cdots\ts m}(u)$ acts on the basis vectors
$\xiL$ by scalar multiplication; see \eqref{agt}. Furthermore, 
$E_{ij}$ with $i<j$ is a commutator in the generators $E_{k,k+1}$
with $k=1,\dots,n-1$.
By the Gelfand--Tsetlin formulas \eqref{Egt}, $E_{k,k+1}\ts\xi^{}_{\La'}$
is a linear combination of the basis vectors
$\xi^{}_{\La'+\delta_{ka}}$, where $a=1,\dots,k$.
The proof is completed by the application of the induction hypothesis to
the pattern $\La'$. 
\endproof

\ble\label{lem:sec}
The $(n-1)$-th row of the pattern $\La^0$
is $(\la_1,\dots,\la_{n-1})$. 
\ele

\Proof We shall be proving
the following property of $\La^0$ which
clearly implies the statement.
For every 
$r=1,\dots,n-1$ we have: if $i\geq r$ and
\beq
\la^0_{n-1,i}= \la^0_{n-2,i-1}=\cdots =	\la^0_{n-r,i-r+1},
\non
\end{equation}
then either $\la^0_{n-1,i}=\la_i$, or $i\geq r+1$ and
\beq
\la^0_{n-1,i}= \la^0_{n-2,i-1}=\cdots =	\la^0_{n-r,i-r+1}=\la^0_{n-r-1,i-r}.
\non
\end{equation}
Suppose the contrary, and let $r$ take the minimum value for which the
property fails. That is, there exists $i\geq r$ such that
$\La':=\La^0+\delta_{n-1,i}+\cdots+\delta_{n-r,i-r+1}$ is a pattern.
Since $\zeta$ is a singular vector, we have
\beql{trz}
{t\ts}^{1\cdots\ts n-r}_{1\cdots\ts n-r-1,n}(u)\ts\zeta=0.
\end{equation}
By Proposition~\ref{prop:delqm},
\beql{trr}
\sum_{c_1<\cdots <c_{n-r}} 
\sum_{\La}{t\ts}^{1\cdots\ts n-r}_{c_1 \ldots \ts c_{n-r}}(u)\ts\xiL
\ot {t\ts}^{c_1\ldots \ts c_{n-r}}_{1\cdots\ts n-r-1,n}(u)\ts\etaL=0.
\end{equation}
The coefficient
of the vector
$\xi^{}_{\La'}\ot\etaLo$
in the expansion of the left hand side of \eqref{trr} must be $0$.
Let us determine which patterns $\La$ yield a nontrivial contribution
to this coefficient. Considering the weight of
${t\ts}^{1\cdots\ts n-r}_{c_1 \ldots \ts c_{n-r}}(u)\ts\xiL$
we come to the relation
\beq
w(\La)+\ve_1+\cdots+\ve_{n-r}-\ve_{c_1}-\cdots-\ve_{c_{n-r}}=
w(\La^0)+\ve_{n-r}-\ve_n,
\non
\end{equation}
and hence
\beq
w(\La)-w(\La^0)=\ve_{c_1}+\cdots+\ve_{c_{n-r}}
-\ve_1-\cdots-\ve_{n-r-1}-\ve_n.
\non
\end{equation}
Since $w(\La)\succeq w(\La^0)$ by Lemma~\ref{lem:onep}, 
we obtain from \eqref{eqp} that
$c_a=a$ for $a=1,\dots,n-r-1$, and $c_{n-r}\in \{n-r,\dots,n\}$.
Then $w(\La)=w(\La^0)+\ve_{c_{n-r}}-\ve_n$. By Lemma~\ref{lem:dif},
$\La$ should be obtained from $\La^0$ by increasing 
exactly one entry by $1$ in each of the rows $c_{n-r}, c_{n-r}+1,\dots,n-1$.
On the other hand, by the minimality of $r$
and the betweenness conditions \eqref{betw}, the array $\La$ would not
be a pattern, unless $c_{n-r}=n-r$ or $c_{n-r}=n$. Thus, the coefficient in question
can only have a contribution from two summands in \eqref{trr}, namely,
\beql{cs1}
{t\ts}^{1\cdots\ts n-r}_{1\cdots\ts n-r-1,n}(u)\ts
\xi^{}_{\La^0}
\ot {t\ts}^{1\cdots\ts n-r-1,n}_{1\cdots\ts n-r-1,n}(u)
\ts\eta^{}_{\La^0}
\end{equation}
and
\beql{cs2}
{t\ts}^{1\cdots\ts n-r}_{1\cdots\ts n-r}(u)\ts
\xi^{}_{\La'}
\ot {t\ts}^{1\cdots\ts n-r}_{1\cdots\ts n-r-1,n}(u)
\ts\eta^{}_{\La'}.
\end{equation}
We consider \eqref{cs1} first. By
\eqref{qmrel},
\beql{comm}
{t\ts}^{1\cdots\ts n-r}_{1\cdots\ts n-r-1,n}(u)=
[{t\ts}^{1\cdots\ts n-r}_{1\cdots\ts n-r}(u),E_{n-r,n}].
\end{equation}
It has been observed above that the minimality of $r$ and the betweenness
conditions \eqref{betw} imply that 
$E_{cn}\ts \xi^{}_{\La^0}=0$ for $n-r< c\leq n-1$.
Hence using the relation
\beq
E_{n-r,n}=[E_{n-r,n-r+1},\cdots[E_{n-3,n-2},
\ts [E_{n-2,n-1},E_{n-1,n}]\ts]\cdots]
\non
\end{equation}
we get
\beql{ta}
E_{n-r,n}\ts\xiLo 
=(-1)^{r-1}E_{n-1,n}E_{n-2,n-1}\cdots E_{n-r,n-r+1} \ts\xiLo.
\non
\end{equation}
Therefore, by \eqref{Egt} the expansion of $E_{n-r,n}\ts\xiLo$
in terms of the basis vectors $\xi^{}_{\La}$ contains
$\xi^{}_{\La'}$ with a nonzero coefficient $C$.
It will now be convenient to use polynomial quantum minor operators
defined by
\beql{Tt}
{T\ts}^{a_1\cdots\ts a_{m}}_{b_1\cdots\ts b_{m}}(u)=
u(u-1)\cdots (u-m+1)\ts
{t\ts}^{a_1\cdots\ts a_{m}}_{b_1\cdots\ts b_{m}}(u),
\non
\end{equation}
see \eqref{qminor}.
Using \eqref{agt} we find from \eqref{comm} that
the coefficient of $\xi^{}_{\La'}$ in the expansion of 
${T\ts}^{1\cdots\ts n-r}_{1\cdots\ts n-r-1,n}(u)\ts\xiLo$ equals
\beq
C\ts (u+l^{\ts 0}_{n-r,1})\cdots\wedge_j\cdots 
(u+l^{\ts 0}_{n-r,n-r}),\qquad j:=i-r+1,
\non
\end{equation}
where $C$ is a nonzero constant.
For the second factor in \eqref{cs1} we find from \eqref{qminor}
and Lemma~\ref{lem:lao}
that
\beql{sfa}
{T\ts}^{1\cdots\ts n-r-1,n}_{1\cdots\ts n-r-1,n}(u)\ts\etaLo
=(u+m_1)\cdots (u+m_{n-r-1})(u+m_n+r)\ts \etaLo.
\non
\end{equation}

Consider now the expression \eqref{cs2}. By \eqref{agt} we have
\beq
{T\ts}^{1\cdots\ts n-r}_{1\cdots\ts n-r}(u)\ts\xi^{}_{\La'}
=(u+l^{\ts 0}_{n-r,1})\cdots(u+l^{\ts 0}_{n-r,j}+1)\cdots 
(u+l^{\ts 0}_{n-r,n-r})\ts \xi^{}_{\La'}.
\non
\end{equation}
Since $\zeta$ is a $\gl_n$-weight vector and $w(\La')=w(\La^0)+\ve_{n-r}-\ve_n$,
the vector $\eta^{}_{\La'}$ is a linear combination
of the $\xi'_{M}$ with
$w(M)=\mu-\ve_{n-r}+\ve_n$,
where $\{\xi'_{M}\}$ is the Gelfand--Tsetlin basis of $L(\mu)$.
This implies that the $(n-r)$-th row of each of the patterns $M$ is
$(\mu_1,\dots,\mu_{n-r-1},\mu_{n-r}-1)$. We therefore have
$E_{n-r,n}\ts \eta^{}_{\La'}=\text{const}\cdot \xi'$
and so, by \eqref{agt} and \eqref{comm}
\beql{tjp}
{T\ts}^{1\cdots\ts n-r}_{1\cdots\ts n-r-1,n}(u)
\ts\eta^{}_{\La'}=
\text{\rm const}\cdot (u+m_1)\cdots (u+m_{n-r-1})\ts \xi'.
\non
\end{equation}

Combining the results of the above calculations, and taking the
coefficient of the vector $\xi^{}_{\La'}\ot\etaLo$ in \eqref{trz},
we obtain
\beq
\bal
C\ts (u+l^{\ts 0}_{n-r,1})\cdots\wedge_j\cdots 
(u+l^{\ts 0}_{n-r,n-r})
(u+m_1)\cdots (u+m_{n-r-1})&(u+m_n+r)\\
{}+\text{const}\cdot(u+l^{\ts 0}_{n-r,1})\cdots(u+l^{\ts 0}_{n-r,j}+1)\cdots 
(u+l^{\ts 0}_{n-r,n-r})&\\
{}\times(u+m_1)\cdots (u+m_{n-r-1})&=0.
\eal
\non
\end{equation}
Deleting the common factors gives
\beq
C\ts (u+m_n+r)
+\text{const}\cdot(u+l^{\ts 0}_{n-r,j}+1)=0.
\non
\end{equation}
Put $u=	-l^{\ts 0}_{n-r,j}-1$ in this relation.  
Since $C$ is nonzero we get
$m_n=l^{\ts 0}_{n-r,j}-r+1$ and thus
$m_n=l^{\ts 0}_{n-1,i}$.
By the betweenness conditions for $\La^0$ and $\La'$,
\beql{betli}
\la_i-\la^0_{n-1,i}>0\quad\text{and}\quad \la^0_{n-1,i}-\la_{i+1}\geq 0,
\non
\end{equation}
which implies that both differences
$l_i-m_n$ and $m_n-l_{i+1}$ are positive integers.
Thus $m_n\in\lan l_{i+1},l_i\ran\subseteq \lan l_{n},l_1\ran$
which contradicts \eqref{1ncond}.
Therefore, our assumption that $\La'$
is a pattern must be wrong. 
\endproof

By Lemma~\ref{lem:dif} we can now conclude that all
vectors $\xiL$ which occur in \eqref{zeta} belong to
the $\U(\gl_{n-1})$-span of the highest vector $\xi$ of $L(\la)$.
This span is isomorphic to the irreducible representation $L(\la_-)$ of
$\gl_{n-1}$ with the highest weight $\la_-=(\la_1,\dots,\la_{n-1})$.
In particular, $E_{nn}\ts\xiL=\la_n\ts\xiL$ for each $\La$.
Furthermore, if $\zeta$ is a linear combination of vectors
$\xi_{\La}\ot \xi'_{M}$
then by Lemma~\ref{lem:lao}, for the corresponding patterns we have
\beql{weilc}
w(\La)+w(M)=w(\La^0)+\mu.
\non
\end{equation}
Therefore, $E_{nn}\ts\xi'_{M}=\mu_n\ts\xi'_{M}$ for all $M$ which
implies that the $(n-1)$-th row of each pattern
$M$ coincides with $(\mu_1,\dots,\mu_{n-1})$.
In other words, each vector $\xi'_{M}$ belongs to
the $\U(\gl_{n-1})$-span of	 $\xi'$ which is isomorphic
to $L(\mu_-)$ where $\mu_-=(\mu_1,\dots,\mu_{n-1})$.
Thus, $\zeta$ belongs to
\beql{zelm}
L(\la_-)\ot L(\mu_-).
\end{equation}
By \eqref{copr} and the defining relations \eqref{defrel}, the $\Y(n-1)$-module structure
on \eqref{zelm} coincides with the one obtained by
restriction from $\Y(n)$ to the subalgebra generated by the $t_{ij}(u)$
with $1\leq i,j\leq n-1$.
The vector $\zeta$ is annihilated by the operators $b_1(u),\dots,b_{n-2}(u)$.
By the assumption of the theorem, for each pair
$(i,j)$ such that $1\leq i<j\leq n-1$
the condition \eqref{ijcond} is satisfied.
Therefore, the $\Y(n-1)$-module \eqref{zelm} 
is irreducible by the induction
hypothesis, and we may finally conclude from Proposition~\ref{prop:sing}
that \eqref{sing} holds.

Next, we derive a similar result for the singular
lowest vectors under the assumptions
of Theorem~\ref{thm:sufc}. As we pointed out,
the condition \eqref{1ncond} can also be assumed due to
Proposition~\ref{prop:perm}.

\ble\label{lem:loww}
If $\zeta'\in L(\la)\ot L(\mu)$ and $t_{ij}(u)\ts\zeta'=0$ for all
$1\leq j<i\leq n$ then
\beql{singl}
\zeta'={\rm const}\cdot\eta\ot\eta',
\end{equation}
where $\eta$ and $\eta'$ are the lowest vectors of $L(\la)$
and $L(\mu)$, respectively.
\ele

\Proof
Let $\omega$ be the permutation of the indices $1,\dots,n$ such that
$\omega(i)=n-i+1$. The mapping
\beql{omeY}
\Y(n)\to \Y(n),\qquad
t_{ij}(u)\mapsto t_{\omega(i)\omega(j)}(u)
\end{equation}
defines an automorphism of the Yangian $\Y(n)$. This follows easily from
the defining relations \eqref{defrel}. We equip the space 
$L=L(\la)\ot L(\mu)$ with
another structure of
$\Y(n)$-module which is obtained by pulling back through
the automorphism \eqref{omeY}. Denote this new representation by $L^{\omega}$.
Similarly, the mapping
\beq
\U(\gl_n)\to \U(\gl_n),\qquad
E_{ij}\mapsto E_{\omega(i)\omega(j)}
\non
\end{equation}
defines an automorphism of $\U(\gl_n)$. Denote by $L(\la)^{\omega}$
the representation of $\U(\gl_n)$ obtained from $L(\la)$ by pulling back through
this automorphism and extend it to $\Y(n)$ using \eqref{epi0}. 
It follows from \eqref{copr}
that the $\Y(n)$-module $L^{\omega}$ is isomorphic to the
tensor product
$L(\la)^{\omega}\ot L(\mu)^{\omega}$.
The weight of the
lowest vector $\eta$
of $L(\la)$ is $(\la_n,\dots,\la_1)$. Therefore
$\eta$, when regarded as an element of $L(\la)^{\omega}$,
is the highest vector of the weight $\la$.
In particular, the $\U(\gl_n)$-module
$L(\la)^{\omega}$ is isomorphic to $L(\la)$.
Now, $\zeta'$ is a singular vector
of the $\Y(n)$-module $L^{\omega}$. By the proved above claim
for the singular vectors, 
$\zeta'$ is, up to a constant factor,
the tensor product of the highest vectors of 
$L(\la)^{\omega}$ and $L(\mu)^{\omega}$,
that is, \eqref{singl} holds.
\endproof

To complete the proof of the theorem, we need to show that the submodule
of $L=L(\la)\ot L(\mu)$ generated by the tensor product
of the highest vectors
$\zeta=\xi\ot\xi'$
coincides with $L$. For this we introduce a $\Y(n)$-module structure
on the space $L^*$ dual to $L$. It follows immediately from the defining
relations \eqref{defrel} that the mapping
\beq
\sigma:\Y(n)\to\Y(n),\qquad t_{ij}(u)\mapsto t_{ij}(-u),
\non
\end{equation}
defines an anti-automorphism of $\Y(n)$. Now, $L^*$ becomes
a $\Y(n)$-module if we set
\beql{defdu}
(yf)(v)=f(\sigma(y)v),\qquad y\in\Y(n),\quad f\in L^*,\quad v\in L.
\end{equation}
Similarly, the dual space $L(\la)^*$ of the $\gl_n$-module $L(\la)$
can be regarded as a $\gl_n$-module with the action defined by
\beq
(E_{ij}f)(v)=f(-E_{ij}v), \qquad f\in L(\la)^*,\quad v\in L(\la).
\non
\end{equation}
We obtain easily from \eqref{copr} that the $\Y(n)$-module $L^*$
is isomorphic to the tensor product
$
L(\la)^*\ot L(\mu)^*,
$
where $L(\la)^*$ and $L(\mu)^*$ are extended to $\Y(n)$
by \eqref{epi0}.
The vector $\xi^{*}\in L(\la)^*$, 
dual to the highest vector
$\xi$, is the lowest vector with the weight 
$-\la$. The highest weight of $L(\la)^*$
will be therefore $-\la^{\omega}=(-\la_n,\dots,-\la_1)$.
Thus, we have
\beql{isolst}
L^*\simeq L(-\la^{\omega})\ot L(-\mu^{\omega}).
\end{equation}

If we assume that the vector $\zeta$ generates a proper submodule $N$ in $L$
then its annihilator
\beq
\text{\rm Ann\ts} N=\{f\in L^*\ |\ f(v)=0\quad\text{\rm for all}\quad v\in N\}
\non
\end{equation}
is a nonzero submodule in $L^*$. 
Hence,  $\text{\rm Ann\ts} N$ must contain
a vector $\zeta'$ which is annihilated by the generators $t_{ij}(u)$ with
$i>j$. However, the condition \eqref{1ncond} remains
satisfied when $\la$ and $\mu$ are respectively
replaced with  $-\la^{\omega}$ and $-\mu^{\omega}$.
So, by Lemma~\ref{lem:loww}, the vector
$\zeta'$ must be, up to a constant factor, the tensor product
of the lowest vectors of the representations $L(-\la^{\omega})$
and $L(-\mu^{\omega})$.
But the vector $\xi^{*}{}\ot\ts {\xi'}^{*}$ does not belong
to $\text{\rm Ann\ts} N$. This makes a contradiction and so,
the submodule generated by $\zeta$ must coincide with $L$. This
completes the proof of Theorem~\ref{thm:sufc}.

\section{Necessary conditions}\label{sec:necc}
\setcounter{equation}{0}

We keep using the notation \eqref{lanran}. As in the previous section,
we assume that $\la$ and $\mu$ are complex $\gl_n$-highest weights.

\bth\label{thm:necc}
Suppose that the $\Y(n)$-module $L(\la)\ot L(\mu)$ is irreducible.
Then for each pair of indices $1\leq i<j\leq n$ we have
\beql{ijcondn}
m_j,m_i\not\in\lan l_j,l_i\ran\qquad\text{or}
\qquad l_j,l_i\not\in\lan m_j,m_i\ran.
\end{equation}
\eth

The proof will follow from a sequence of lemmas. We use induction
on $n$. 
Given a $\gl_n$-highest weight $\la=(\la_1,\dots,\la_n)$ we set
\beql{lapm}
\la_-=(\la_1,\dots,\la_{n-1})\quad\text{and}\quad\la_+=(\la_2,\dots,\la_{n}).
\non
\end{equation}

\ble\label{lem:lapm}
If the $\Y(n)$-module $L(\la)\ot L(\mu)$ is irreducible
then both $\Y(n-1)$-modules $L(\la_-)\ot L(\mu_-)$ and $L(\la_+)\ot L(\mu_+)$
are irreducible.
\ele

\Proof We shall identify $L(\la_-)$ and $L(\mu_-)$ with
the $\U(\gl_{n-1})$-spans of the highest vectors $\xi$ in $L(\la)$
and $\xi'$ in $L(\mu)$, respectively.
Any generator $E_{in}$ of $\gl_n$ with $i<n$ annihilates $L(\la_-)$ and $L(\mu_-)$.
Hence, by \eqref{copr} and \eqref{epi0}, 
the subspace $L(\la_-)\ot L(\mu_-)$ of $L(\la)\ot L(\mu)$
is invariant with respect to the action of the subalgebra $\Y(n-1)$ of $\Y(n)$,
and this action coincides with the one defined in Section~\ref{sec:pre}.

Suppose that there is a nonzero submodule in
$L(\la_-)\ot L(\mu_-)$
which does not contain the vector $\xi\ot\xi'$.  Then this submodule contains
a $\Y(n-1)$-singular vector $\zeta$. However, $\zeta$ must also be a 
$\Y(n)$-singular vector. Indeed, $t_{in}(u)\ts\zeta=0$
for any $i<n$ which easily follows from \eqref{defrel} and \eqref{copr}.
This implies that $L(\la)\ot L(\mu)$ is not irreducible,
contradiction.

Suppose now that the $\Y(n-1)$-submodule of
$L_-=L(\la_-)\ot L(\mu_-)$ generated by $\xi\ot\xi'$ is proper.
It follows from the defining
relations \eqref{defrel} that the mapping
\beql{autr}
\tau:\Y(n-1)\to\Y(n-1),\qquad t_{ij}(u)\mapsto t_{ji}(u),
\non
\end{equation}
defines an anti-automorphism of $\Y(n-1)$. The dual space 
${\displaystyle L^{\ts *}_-}$ becomes
a $\Y(n-1)$-module if we set
\beql{defdum}
(yf)(v)=f(\tau(y)v),\qquad y\in\Y(n-1),\quad f\in L^{\ts *}_-,\quad v\in L_-.
\non
\end{equation}
We easily derive from \eqref{copr} that the module
${\displaystyle L^{\ts *}_-}$ is isomorphic to $L(\mu_-)\ot L(\la_-)$.
Since $\xi\ot\xi'$ generates a proper submodule in $L_-$,
its annihilator in ${\displaystyle L^{\ts *}_-}$ is a nonzero submodule
which does not contain the vector $\xi'\ot \xi\in L(\mu_-)\ot L(\la_-)$.
However, the $\Y(n)$-module $L(\mu)\ot L(\la)$ is irreducible
by Proposition~\ref{prop:perm}. Thus, our assumption leads again to a contradiction
due to the previous argument. This proves that the $\Y(n-1)$-module 
$L(\la_-)\ot L(\mu_-)$ is irreducible.

Now consider the $\Y(n-1)$-module $L(\la_+)\ot L(\mu_+)$.
If the $\Y(n)$-module $L=L(\la)\ot L(\mu)$ is irreducible then so is
the module $L^*$ defined in \eqref{defdu}. To complete the proof
we apply the isomorphism \eqref{isolst} and
the above argument. \endproof

Due to Lemma~\ref{lem:lapm},
if the $\Y(n)$-module $L(\la)\ot L(\mu)$ is irreducible then,
by the induction hypothesis, both conditions 
$\Cond(\A_{\la}^-,\A_{\mu}^-)$ and $\Cond(\A_{\la}^+,\A_{\mu}^+)$
hold. Therefore, by Proposition~\ref{prop:equiv} the
conditions
\eqref{ijcondn} are satisfied
for all pairs $(i,j)\ne (1,n)$. Suppose that they 
are violated for the pair $(1,n)$. 
Then, as was shown in the proof of Proposition~\ref{prop:equiv},
the condition \eqref{poss} should hold. Using Proposition~\ref{prop:perm}, 
if necessary, we may assume that
$m_n\in\lan l_n,l_1\ran$ and
$l_1\in\lan m_n,m_1 \ran$.
Therefore, there exist indices $p,q\in \{1,\dots,n-1\}$ such that
\beql{pcond}
m_n\in\lan l_{p+1},l_p\ran\qquad\text{and}\qquad l_1\in\lan m_{q+1},m_q\ran.
\end{equation}
If $p=n-1$ then by
$\Cond(\A_{\la}^+,\A_{\mu}^+)$ we must have
$q=1$; cf. the proof of Proposition~\ref{prop:equiv}.
Thus, $l_1\in\lan m_{2},m_1\ran$.
Moreover, using also the condition $\Cond(\A_{\la}^-,\A_{\mu}^-)$,
we conclude that there should exist indices $r$ and $s$ such that
\beql{indrs}
m_2,\dots,m_r\in\{l_2,\dots,l_s\},\qquad
l_{s+1},\dots,l_{n-1}\in\{m_{r+1},\dots,m_{n-1}\},
\non
\end{equation}
as shown in the picture:

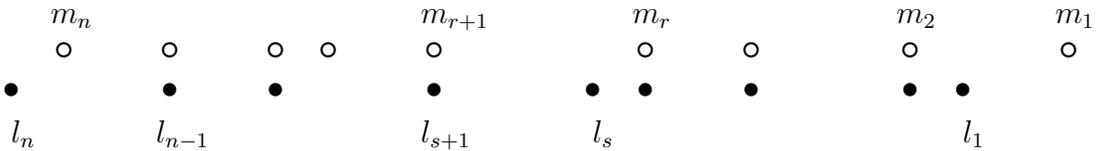
\begin{center}
\begin{picture}(400,60)
\thicklines
\put(0,20){\circle*{5}}
\put(60,20){\circle*{5}}
\put(100,20){\circle*{5}}
\put(160,20){\circle*{5}}
\put(220,20){\circle*{5}}
\put(240,20){\circle*{5}}
\put(280,20){\circle*{5}}
\put(340,20){\circle*{5}}
\put(360,20){\circle*{5}}

\put(0,0){$l_n$}
\put(155,0){$l_{s+1}$}
\put(220,0){$l_{s}$}
\put(55,0){$l_{n-1}$}
\put(360,0){$l_{1}$}

\put(20,35){\circle{5}}
\put(60,35){\circle{5}}
\put(100,35){\circle{5}}
\put(120,35){\circle{5}}
\put(160,35){\circle{5}}
\put(240,35){\circle{5}}
\put(280,35){\circle{5}}
\put(340,35){\circle{5}}
\put(400,35){\circle{5}}

\put(15,45){$m_n$}
\put(155,45){$m_{r+1}$}
\put(235,45){$m_{r}$}
\put(335,45){$m_{2}$}
\put(395,45){$m_1$}

\end{picture}
\end{center}

\noindent
In particular, this implies that
\beql{ineqml}
l_i-m_i\in\ZZ_+ \qquad\text{for all}\quad i=2,\dots,n-1.
\end{equation}
Now let $p\leq n-2$ in \eqref{pcond}.
The conditions $\Cond(\A_{\la}^-,\A_{\mu}^-)$ and $\Cond(\A_{\la}^+,\A_{\mu}^+)$
imply that
\beql{leqm}
l_{p-i+1}=m_{n-i}\qquad\text{for}\quad i=1,\dots,p-1
\end{equation}
while
\beql{lone}
l_1\in\lan m_{n-p+1},m_{n-p}\ran,
\end{equation}
as illustrated:

\begin{center}
\begin{picture}(400,60)
\thicklines
\put(0,20){\circle*{5}}
\put(40,20){\circle*{5}}
\put(60,20){\circle*{5}}
\put(100,20){\circle*{5}}
\put(160,20){\circle*{5}}
\put(200,20){\circle*{5}}
\put(220,20){\circle*{5}}
\put(240,20){\circle*{5}}
\put(280,20){\circle*{5}}

\put(0,0){$l_n$}
\put(100,0){$l_{p+1}$}
\put(160,0){$l_{p}$}
\put(240,0){$l_{2}$}
\put(280,0){$l_{1}$}

\put(120,35){\circle{5}}
\put(160,35){\circle{5}}
\put(200,35){\circle{5}}
\put(220,35){\circle{5}}
\put(240,35){\circle{5}}
\put(300,35){\circle{5}}
\put(340,35){\circle{5}}
\put(380,35){\circle{5}}
\put(400,35){\circle{5}}

\put(110,45){$m_n$}
\put(150,45){$m_{n-1}$}
\put(230,45){$m_{n-p+1}$}
\put(290,45){$m_{n-p}$}
\put(395,45){$m_1$}

\end{picture}
\end{center}

Let $L$ be a highest weight module over $\Y(n)$
generated by a vector $\zeta$ such that
\beql{tiiarhw}
T_{ii}(u)\ts\zeta=(u+\la_i)(u+\mu_i)\ts\zeta,\qquad i=1,\dots,n,
\end{equation}
where $T_{ij}(u)=u^2\ts t_{ij}(u)$. We shall also suppose that
the elements $t_{ij}^{(r)}$ with $r\geq 3$ act trivially on $L$
so that the operators $T_{ij}(u)$ are polynomials in $u$.
In other words, $L$ is a module over the 
quotient algebra $\Y(n)/I$ where $I$ is the ideal generated by the elements 
$t_{ij}^{(r)}$ with $r\geq 3$. 

Given sequences $a_1,\dots, a_m$ and $b_1,\dots,b_m$ 
of elements of $\{1,\dots,n\}$ we denote by 
${T\ts}^{a_1\cdots\ts a_{m}}_{b_1\cdots\ts b_{m}}(u)$
the corresponding quantum minor defined by
\eqref{qminor} or \eqref{qminor2} with $t_{ij}(u)$ 
respectively replaced by $T_{ij}(u)$.
Similarly, we define the operators $B_m(u)$ and $A_m(u)$ 
by \eqref{drgen} with the same replacement.

For $1\leq a< r\leq n$ introduce the 
{\it raising operators\/} $\tau_{ar}(v)$ and
{\it lowering operators\/}
$\tau_{ra}(v)$ \cite{m:gt} by
\beql{tau}
\tau^{}_{ar}(v)={T\ts}^{1\ts\cdots\ts a}_{1\ts\cdots\ts a-1,r}(v),
\qquad
\tau^{}_{ra}(v)={T\ts}^{a+1\ts\cdots\ts r}_{a\ts\cdots\ts r-1}(v),
\non
\end{equation}
where $v$ is a variable. We also set $\tau^{}_{ra}(v)\equiv 1$ for $r\leq a$.
For any non-negative integer $k$ introduce the product of the
lowering operators
\beql{Tau}
\Tau^{}_{ra}(v,k)=\tau^{}_{ra}(v+k-1)\cdots\tau^{}_{ra}(v+1)
\ts\tau^{}_{ra}(v).
\end{equation}
Suppose now that $\eta\in L$ is a $\Y(n-1)$-singular vector. That is,
$\eta$ is annihilated by $T_{ij}(u)$ for $1\leq i<j\leq n-1$ and
\beql{tiieta1}
T_{ii}(u)\ts\eta=\nu_i(u)\ts\eta,\qquad i=1,\dots,n-1,
\end{equation}
for some polynomials $\nu_i(u)$ of degree two.
The next three lemmas will provide a basis for our calculations.

\ble\label{lem:tii+1} 
We have the following relations in $L$:
\beql{tiinaze}
T_{ii}(u)\ts\Tau^{}_{na}(v,k)\ts\eta
=\nu_i(u)\ts\Tau^{}_{na}(v,k)\ts\eta,
\end{equation}
if $1\leq i\leq n-1$ and $i\ne a$, while
\begin{align}\label{taanaze}
T_{aa}(u)\ts \Tau^{}_{na}(v,k)\ts\eta=\frac{(u-v-k)\ts\nu_a(u)}{u-v}
\ts  &\Tau^{}_{na}(v,k)\ts\eta\\
\non
{}+\frac{k}{u-v}\ts 
\sum_{c=a+1}^n\nu_a(v)&\ts\nu_{a+1}(v-1)\cdots\nu_{c-1}(v-c+a+1)\\
\non
{}\times{}&
T_{ca}(u)\ts \Tau^{}_{na}(v+1,k-1)\ts\tau_{nc}(v-c+a)\ts\eta.
\end{align}
Moreover,
\beql{tii+1na}
T_{i,i+1}(u)\ts\Tau^{}_{na}(v,k)\ts\eta=0,
\end{equation}
if $1\leq i<n-1$ and $i\ne a$, while for $a<n-1$
\beql{taa+1naze}
\bal
T_{a,a+1}(u)\ts \Tau^{}_{na}(v,k)\ts\eta=\frac{k}{u-v}\ts 
\sum_{c=a+1}^n\nu_a(v)&\ts\nu_{a+1}(v-1)\cdots\nu_{c-1}(v-c+a+1)\\
{}\times{}&
T_{c,a+1}(u)\ts \Tau^{}_{na}(v+1,k-1)\ts\tau^{}_{nc}(v-c+a)\ts\eta.
\eal  
\end{equation}
In particular, if $\nu_a(-\rho)=0$ for some $\rho$ then
\beql{taaa}
T_{aa}(u)\ts\Tau^{}_{na}(-\rho,k)\ts\eta
=\frac{(u+\rho-k)\ts\nu_a(u)}{u+\rho}\ts\Tau^{}_{na}(-\rho,k)\ts\eta,
\end{equation}
and for $a<n-1$
\beql{taaa+1}
T_{a,a+1}(u)\ts\Tau^{}_{na}(-\rho,k)\ts\eta=0.
\end{equation}
\ele

\Proof Note that the coefficients of $\Tau^{}_{na}(v,k)$ 
are linear combinations of monomials in the generators $t_{rs}^{(j)}$ with
$a\leq s\leq n-1$. 
Suppose that $i<a$. We have $T_{il}(u)\ts\eta=0$ for $a\leq l\leq n-1$.
Therefore, applying \eqref{defrel} to the commutators $[T_{ii}(u),T_{rs}(v)]$
and $[T_{il}(u),T_{rs}(v)]$,
we conclude by an easy induction 
that $T_{il}(u)\ts \Tau^{}_{na}(v,k)\ts\eta=0$.
This proves \eqref{tiinaze} and
\eqref{tii+1na} for $i<a$.	For $i>a$ both these relations
are immediate from \eqref{center}. Further, by \eqref{qmrel},
\beql{tpptau}
\bal
{T}^{}_{aa}(u)\ts\Tau^{}_{na}(v,k)
{}=\frac{u-v-k}
{u-v-k+1} \ts\tau^{}_{na}(v+k-1)\ts{T}^{}_{aa}(u)\ts&\Tau^{}_{na}(v,k-1) \\
{}+\frac{1}
{u-v-k+1} \ts\sum_{c=a+1}^n
{T}^{}_{ca}(u)\ts
{T\ts}^{a\ts\cdots\ts \wh c\ts\cdots\ts n}_{a\ts\cdots\ts n-1}(v+k-1)
\ts&\Tau^{}_{na}(v,k-1)
\ts(-1)^{c-a-1}.
\eal
\non
\end{equation} 
The subalgebra $\Y_a$ of $\Y(n)$ 
generated by $t_{rs}(u)$ with $a\leq r,s\leq n$
is naturally isomorphic to the Yangian $\Y(n-a+1)$. 
Applying the automorphism \eqref{audual}
to this subalgebra,
we derive from the defining relations \eqref{defrel}
that
\beql{permin}
{T\ts}^{a\ts\cdots\ts \wh c\ts\cdots\ts n}_{a\ts\cdots\ts n-1}(v+1)
\ts{T\ts}^{a+1\ts\cdots\ts n}_{a\ts\cdots\ts n-1}(v)
={T\ts}^{a+1\ts\cdots\ts n}_{a\ts\cdots\ts n-1}(v+1)
\ts{T\ts}^{a\ts\cdots\ts \wh c\ts\cdots\ts n}_{a\ts\cdots\ts n-1}(v)
\non
\end{equation}
for every $c=a+1,\dots,n$. Hence,
\beql{traintau}
{T\ts}^{a\ts\cdots\ts \wh c\ts\cdots\ts n}_{a\ts\cdots\ts n-1}(v+k-1)
\ts\Tau^{}_{na}(v,k-1)
=\Tau^{}_{na}(v+1,k-1)\ts
{T\ts}^{a\ts\cdots\ts \wh c\ts\cdots\ts n}_{a\ts\cdots\ts n-1}(v).
\non
\end{equation}
Note that
${T}^{}_{ca}(u)$ commutes with $\tau^{}_{na}(v+k-1)$ by \eqref{center}. Therefore,
an easy induction on $k$ gives
\beql{tpptauk}
\bal
{T}^{}_{aa}(u)\ts\Tau^{}_{na}(v,k)
{}&=\frac{u-v-k}
{u-v} \ts\Tau^{}_{na}(v,k)\ts{T}^{}_{aa}(u)\\
{}+\frac{k}
{u-v} \ts&\sum_{c=a+1}^n
{T}^{}_{ca}(u)\ts
\Tau^{}_{na}(v+1,k-1)\ts
{T\ts}^{a\ts\cdots\ts \wh c\ts\cdots\ts n}_{a\ts\cdots\ts n-1}(v)
\ts(-1)^{c-a-1}.
\eal
\end{equation}
The same argument proves the following counterpart of 
\eqref{tpptauk}: for $a<n-1$
\beql{tpptp+1}
\bal
{T}^{}_{a,a+1}(u)\ts\Tau^{}_{na}(v,k)
{}&=\frac{u-v-k}
{u-v} \ts\Tau^{}_{na}(v,k)\ts{T}^{}_{a,a+1}(u)\\
{}+\frac{k}
{u-v} \ts&\sum_{c=a+1}^n
{T}^{}_{c,a+1}(u)\ts
\Tau^{}_{na}(v+1,k-1)\ts
{T\ts}^{a\ts\cdots\ts \wh c\ts\cdots\ts n}_{a\ts\cdots\ts n-1}(v)
\ts(-1)^{c-a-1}.
\eal
\end{equation}
By \eqref{qminor2},
\beql{Tagen}
\bal
{T\ts}^{a\ts\cdots\ts \wh c\ts\cdots\ts n}_{a\ts\cdots\ts n-1}(v)
=\sum_{\sigma\in\Sym_{n-a}}\sgn\sigma&\cdot
T_{n,\sigma(n-1)}(v-n+a+1)\cdots T_{c+1,\sigma(c)}(v-c+a)\\
{}&\times T_{c-1,\sigma(c-1)}(v-c+a+1)\cdots T_{a,\sigma(a)}(v).
\eal
\non
\end{equation} 
Since
$T_{ij}(u)\ts\eta=0$ for $1\leq i<j\leq n-1$ 
we conclude from \eqref{tiieta1} that
\beql{qmca}
{T\ts}^{a\ts\cdots\ts \wh c\ts\cdots\ts n}_{a\ts\cdots\ts n-1}(v)\ts\eta
=\tau^{}_{nc}(v-c+a)\ts\nu_a(v)
\ts\nu_{a+1}(v-1)\cdots\nu_{c-1}(v-c+a+1)\ts\eta.
\non
\end{equation}
The proof is completed by using \eqref{tpptauk} and
\eqref{tpptp+1}.
\endproof

\ble\label{lem:ecomm}
Let $1\leq a<n-1$. Then
we have the relations in $\Y(n)$:
\beql{ecomta}
[E_{n-1,n}\ts\Tau^{}_{na}(v,k)]=
-k \ts\Tau^{}_{na}(v,k-1)\ts 
{T\ts}^{a+1\ts\cdots\ts n}_{a\ts\cdots\ts n-2,n}(v+k-1).
\end{equation}
Moreover,
\beql{tautau}
{T\ts}^{a+1\ts\cdots\ts n}_{a\ts\cdots\ts n-2,n}(u)\ts
{T\ts}^{a+1\ts\cdots\ts n-1}_{a+1\ts\cdots\ts n-1}(u)=
\tau^{}_{n-1,a}(u)\ts{T\ts}^{a+1\ts\cdots\ts n}_{a+1\ts\cdots\ts n}(u)
+ \tau^{}_{na}(u)\ts{T\ts}^{a+1\ts\cdots\ts n-1}_{a+1\ts\cdots\ts n-2,n}(u).
\end{equation}
\ele

\Proof By \eqref{qmrel},
$[E_{n-1,n},\tau^{}_{na}(v)]=-
{T\ts}^{a+1\ts\cdots\ts n}_{a\ts\cdots\ts n-2,n}(v)$.
Since the elements $\tau^{}_{na}(u)$ and $\tau^{}_{na}(v)$ commute,
we can write
\beql{t1Tau0}
[E_{n-1,n}\ts\Tau^{}_{na}(v,k)]
=-\sum_{i=1}^k \tau^{}_{na}(v)\cdots 
{T\ts}^{a+1\ts\cdots\ts n}_{a\ts\cdots\ts n-2,n}(v+i-1)
\cdots \tau^{}_{na}(v+k-1).
\end{equation}
Applying the automorphism  \eqref{audual}
to the subalgebra $\Y_a$
introduced in the proof of Lemma~\ref{lem:tii+1},
we derive from \eqref{defrel}
that
\beql{permtt0}
{T\ts}^{a+1\ts\cdots\ts n}_{a\ts\cdots\ts n-2,n}(v)\ts\tau^{}_{na}(v+1)
=\tau^{}_{na}(v)\ts {T\ts}^{a+1\ts\cdots\ts n}_{a\ts\cdots\ts n-2,n}(v+1).
\end{equation}
Together with \eqref{t1Tau0} this proves \eqref{ecomta}
by an easy induction.

To prove \eqref{tautau} consider the expression
provided by Proposition~\ref{prop:qmr} for the commutator
$[\tau^{}_{n-1,a}(u), {T\ts}^{a+1\ts\cdots\ts n}_{a+1\ts\cdots\ts n}(v)]$.
Multiply both sides of the relation by $u-v$ and put $u=v$.
The terms which do not vanish after this operation correspond to
the maximum value of the summation parameter in the formula,
which implies \eqref{tautau}.
\endproof

We keep using the notation of Lemma~\ref{lem:tii+1}. As before,
$\zeta$ is the highest vector of $L$ satisfying
\eqref{tiiarhw}. We also regard $L$ as a $\gl_n$-module
using \eqref{emb}.

\ble\label{lem:yn-1} We have
the following relations in $L$. 
\beql{ynknle}
E_{in}\ts\Tau^{}_{na}(v,k)\ts\zeta=0\qquad \text{if}\quad i<a.
\end{equation}
If $a<i\leq n-1$ then
\beql{ynknnne}
E_{in}\ts\Tau^{}_{na}(v,k)\ts\zeta
=(-1)^{n-i}\ts k \ts\prod_{j=i+1}^n(v'+l_j)(v'+m_j)
\ts\tau^{}_{ia}(v+k-1)\ts \Tau^{}_{na}(v,k-1)\ts\zeta,
\end{equation}
where $v'=v+a+k-1$. Moreover, if $1\leq a\leq n-1$ then
\begin{multline}\label{ynknn-1}
E_{an}\ts\Tau^{}_{na}(v,k)\ts\zeta
=(-1)^{n-a-1}\ts k\\
{}\times\Big( \Tau^{}_{na}(v+1,k-1)\ts
{T\ts}^{a\ts\cdots\ts n-1}_{a\ts\cdots\ts n-1}(v)
-\Tau^{}_{na}(v,k-1)\ts
{T\ts}^{a+1\ts\cdots\ts n}_{a+1\ts\cdots\ts n}(v+k-1)\Big)\ts\zeta.
\end{multline}
\ele

\Proof The relation \eqref{ynknle} follows from the fact that
the coefficients of $\Tau^{}_{na}(v,k)$ belong to the subalgebra $\Y_a$:
see the proof of Lemma~\ref{lem:tii+1}.
We easily deduce from \eqref{qmrel}
that
\beql{t1tau}
[E_{an},\tau^{}_{na}(v)]=(-1)^{n-a-1}\ts
{T\ts}^{a\ts\cdots\ts n-1}_{a\ts\cdots\ts n-1}(v)
-(-1)^{n-a-1}\ts{T\ts}^{a+1\ts\cdots\ts n}_{a+1\ts\cdots\ts n}(v).
\non
\end{equation}
Therefore,
\begin{multline}\label{t1Tau}
E_{an}\ts\Tau^{}_{na}(v,k)\ts\zeta
=(-1)^{n-a-1}\\
{}\times\sum_{i=1}^k \tau^{}_{na}(v)\cdots 
\big({T\ts}^{a\ts\cdots\ts n-1}_{a\ts\cdots\ts n-1}\ts(v+i-1) 
-{T\ts}^{a+1\ts\cdots\ts n}_{a+1\ts\cdots\ts n}\ts(v+i-1)\big)
\cdots \tau^{}_{na}(v+k-1)\ts\zeta.
\end{multline}
Furthermore, by analogy with \eqref{permtt0} we get
\beql{permtt}
{T\ts}^{a+1\ts\cdots\ts n}_{a+1\ts\cdots\ts n}\ts(v)\ts\tau^{}_{na}(v+1)
=\tau^{}_{na}(v)\ts {T\ts}^{a+1\ts\cdots\ts n}_{a+1\ts\cdots\ts n}\ts(v+1).
\non
\end{equation}
The expression \eqref{t1Tau} now takes the form
\begin{multline}
(-1)^{n-a-1}
\ts\sum_{i=1}^k \tau^{}_{na}(v)\cdots \tau^{}_{na}(v+i-2)\ts
{T\ts}^{a\ts\cdots\ts n-1}_{a\ts\cdots\ts n-1}\ts(v+i-1) 
\ts\Tau^{}_{na}(v+i,k-i)\ts\zeta\\
-(-1)^{n-a-1}\ts k\ts \Tau^{}_{na}(v,k-1)\ts
{T\ts}^{a+1\ts\cdots\ts n}_{a+1\ts\cdots\ts n}\ts(v+k-1)\ts\zeta.
\non
\end{multline}
Similarly, applying again \eqref{audual}
to $\Y_a$ we bring the sum here by an easy induction to the form
\begin{multline}
\sum_{i=1}^k\Big((k-i+1)\ts \tau^{}_{na}(v+k-1)\cdots\wh{\tau^{}_{na}(v+i-1)}
\cdots \tau^{}_{na}(v)\ts {T\ts}^{a\ts\cdots\ts n-1}_{a\ts\cdots\ts n-1}\ts(v+i-1)\\
-(k-i)\ts \tau^{}_{na}(v+k-1)\cdots\wh{\tau^{}_{na}(v+i)}
\cdots \tau^{}_{na}(v)\ts {T\ts}^{a\ts\cdots\ts n-1}_{a\ts\cdots\ts n-1}\ts(v+i)
\Big)\ts\zeta,
\non
\end{multline}
which simplifies to
$k\ts \Tau^{}_{na}(v+1,k-1)\ts
{T\ts}^{a\ts\cdots\ts n-1}_{a\ts\cdots\ts n-1}\ts(v)\ts\zeta$
thus proving \eqref{ynknn-1}.

By \eqref{qmrel}, if $i>a$ then
$[E_{in},\tau^{}_{na}(v)]=(-1)^{n-i}\ts 
{T\ts}^{a+1\ts\cdots\ts n}_{a\ts\cdots\ts\wh{i}\ts\cdots\ts n}(v)$.
As in the proof of \eqref{ecomta}, this brings
the left hand side of \eqref{ynknnne} to the form
\beql{summule}
(-1)^{n-i}\ts k\ts \Tau^{}_{na}(v,k-1)\ts
{T\ts}^{a+1\ts\cdots\ts n}_{a\ts\cdots\ts\wh{i}\ts\cdots\ts n}(v+k-1)\ts\zeta.
\non
\end{equation}
Using \eqref{qminor}, we get
\beql{summex}
{T\ts}^{a+1\ts\cdots\ts n}_{a\ts\cdots\ts\wh{i}\ts\cdots\ts n}(v+k-1)\ts\zeta
=\prod_{j=i+1}^n(v'+l_j)(v'+m_j)
\ts\tau^{}_{ia}(v+k-1)\ts\zeta,
\non
\end{equation}
which completes the proof.
\endproof

We now consider the irreducible highest weight module $V(\la,\mu)$	
over $\Y(n)$ generated by the highest vector $\zeta$ satisfying
\eqref{tiiarhw}.
It follows easily from \eqref{defrel} and the irreducibility of 
$V(\la,\mu)$ that all elements $t_{ij}^{(r)}$ with $r\geq 3$
act trivially in this module.
Furthermore, $V(\la,\mu)$ is isomorphic to the irreducible quotient
of the submodule of $L(\la)\ot L(\mu)$ generated by the vector $\xi\ot\xi'$;
see \eqref{copr} and \eqref{epi0}.
With the parameter $p$ defined in \eqref{pcond}, 
the numbers
\beql{ki}
k_i=l_i-m_{n-p+i},\qquad i=1,\dots,p
\non
\end{equation}
are positive integers by \eqref{ineqml} and \eqref{leqm}.
Introduce the vector $\theta\in V(\la,\mu)$ by
\beql{singtau}
\theta=\Tau^{}_{n-p+1,1}(-\la_1,k_1)\ts
\Tau^{\ts\prime}_{n-p+2,2}(-\la_2,k_2)\cdots 
\Tau^{\ts\prime}_{np}(-\la_p,k_p)\ts\zeta,
\end{equation}
where $\zeta$ is the highest vector of $V(\la,\mu)$,
and $\Tau^{\ts\prime}_{ra}(v,k_a)$ denotes the derivative of the polynomial 
$\Tau^{}_{ra}(v,k_a)$; see \eqref{Tau}.
Our aim is to prove that the vector $\theta$ is zero.
We shall do this in Lemma~\ref{lem:theta=0} below.
The idea of the proof is to show that $\theta$
is annihilated by the operators
$B_i(u)$ for all $i=1,\dots,n-1$ and then apply Proposition~\ref{prop:sing}
noting that $\theta$ is obviously not proportional to $\zeta$.
This works directly in the case $p=1$.
However, if $p\geq 2$ then applying the operators $B_i(u)$ to $\theta$,
we come to a more general problem to prove that
all vectors parametrized by
a certain finite family of pattern-like arrays $\La$ associated with $\la$
are zero. We prove a preliminary lemma first which
describes the properties of these vectors.
The arrays $\La$ which arise in this way will be called {\it admissible\/}.
They are defined as follows.
Each $\La$ is a sequence of rows $\La_r=(\la_{r1},\dots,\la_{rr})$
with $r=1,\dots,n$
of the form described in Section~\ref{sec:pre}.
The top row $\La_n$ coincides with $\la$ and for all $r$ 
the following conditions hold
\beql{betwha}
\la_{ri}-\la_{r-1,i}\in\ZZ_+\qquad{\rm for}\quad i=1,\dots,r-1.
\non
\end{equation}
Each entry $\la_{ri}$ of $\La$ is equal to $\la_i$ 
unless
\beql{ineqki}
i=2,\dots,p\qquad\text{and}\qquad r<n-p+i.
\end{equation}
Moreover, we also require that if $i\in\{2,\dots,p\}$ then
\beql{eqlala}
l_{ii}-m_{n-p+i}\in\ZZ_+,
\end{equation}
where we denote $l_{ri}=\la_{ri}-i+1$.
This condition implies that $0\leq \la_{ri}-\la_{r-1,i}\leq k_i$
for all $i$.
By definition, only a part
of an admissible array can vary with the remaining entries fixed,
as illustrated:

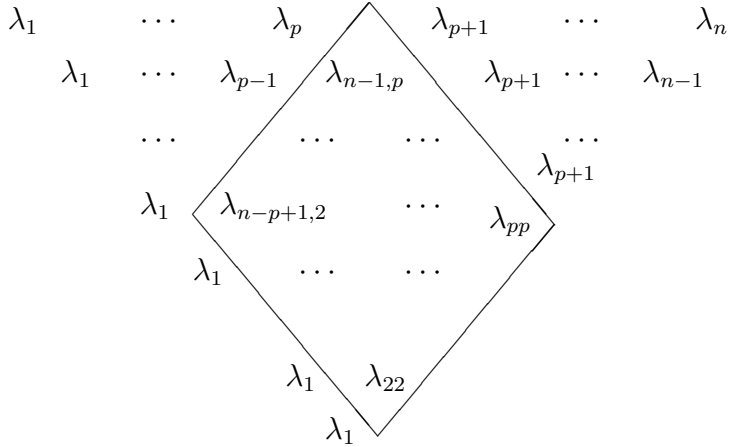
\begin{center}
\begin{picture}(400,200)

\put(60,180){$\la_1$}
\put(110,180){$\cdots$}
\put(160,180){$\la_p$}
\put(220,180){$\la_{p+1}$}
\put(270,180){$\cdots$}
\put(320,180){$\la_n$}

\put(80,160){$\la_1$}
\put(110,160){$\cdots$}
\put(140,160){$\la_{p-1}$}
\put(180,160){$\la_{n-1,p}$}
\put(240,160){$\la_{p+1}$}
\put(270,160){$\cdots$}
\put(300,160){$\la_{n-1}$}

\put(110,135){$\cdots$}
\put(170,135){$\cdots$}
\put(210,135){$\cdots$}
\put(270,135){$\cdots$}

\put(110,110){$\la_1$}
\put(140,110){$\la_{n-p+1,2}$}
\put(210,110){$\cdots$}
\put(260,125){$\la_{p+1}$}

\put(130,85){$\la_1$}
\put(170,85){$\cdots$}
\put(210,85){$\cdots$}

\put(242,105){$\la_{pp}$}


\put(165,45){$\la_1$}
\put(195,45){$\la_{22}$}

\put(180,25){$\la_1$}

\put(130,110){\line(5,-6){70}}
\put(130,110){\line(5,6){66.8}}
\put(200,26){\line(5,6){66.8}}
\put(197,190){\line(5,-6){70}}

\end{picture}
\end{center}

\vspace{-25pt}

\noindent
Given such an array, we set
\beql{thela}
\theta_{\La}=\prod_{r=3,\dots,n\ }^{\rightarrow}\left( \ts
\prod_{i=2}^{p}\wt{\Tau}^{}_{ri}(-\la_{ri}, \la_{ri}-\la_{r-1,i})\right)\ts\zeta,
\non
\end{equation}
where the polynomials $\wt{\Tau}^{}_{ri}(v,k)$ are defined by
\beql{taitidef}
\wt{\Tau}^{}_{ri}(v,k)=
\begin{cases}
\Tau^{\ts\prime}_{ri}(v,k)\qquad&\text{if}\ \ 
r=n-p+i \ \ \text{and}\ \ k=k_i	\\
\Tau^{}_{ri}(v,k)\qquad&\text{otherwise.}
\end{cases}
\end{equation}
Recall also that $\Tau^{}_{ri}(v,k)\equiv 1$
if $r\leq i$.
Note that the factors in the brackets commute for any index $r$ due to
\eqref{center}. We have
$\zeta=\theta^{}_{\La^0}$
for the array $\La^0$ with $\la^0_{ri}=\la_i$ for all $r$ and $i$. Furthermore,
$\theta=\Tau^{}_{n-p+1,1}(-\la_1,k_1)\ts\theta^{}_{\La}$ for the array $\La$
with $l_{ri}=m_{n-p+i}$ for all indices $r$ and $i$ satisfying \eqref{ineqki}.
We define the {\it weight\/} $w(\La)$ of 
an admissible array $\La$ by \eqref{weigt}. We use 
the ordering on the weights described in Section~\ref{sec:sufc}.
Given $\La$,
take the minimum index $r=r(\La)$ such that 
for some $2\leq a\leq p$ the difference $\la_{ra}-\la_{r-1,a}$
is a positive integer.
The following relations for	the admissible arrays $\La$
with $r(\La)\geq n-p+2$ will be used in Lemma~\ref{lem:theta=0}.

\ble\label{lem:theLa} We have the following relations in $V(\la,\mu)$: 
\begin{alignat}{2}\label{biuthr}
T_{ii}(u)\ts\theta_{\La}&=(u+\la_{r-1,i})(u+\mu_i)\ts\theta_{\La}
\qquad&&\text{for}\quad 1\leq i\leq r-1, \\
\label{biu+1r}
T_{ij}(u)\ts\theta_{\La}&=0
\qquad&&\text{for}\quad 1\leq i<j\leq r-1,\\
\label{bathL2r}
B_i(u)\ts\theta_{\La}&=
\sum_{j=1}^i\beta_{ij}(u,\La)\ts\theta_{\La+\delta^{}_{ij}}
\qquad&&\text{for}\quad r-1\leq i< n,
\end{alignat}
where $\beta_{ij}(u,\La)$
are some polynomials in $u$, and we suppose that $\theta_{\La}=0$
if $\La$ is not an admissible array.
\ele

\Proof  
Let $2\leq a\leq p$ be the least index such that 
$k=\la_{ra}-\la_{r-1,a}>0$. Then
\beql{thlade}
\theta_{\La}=\wt{\Tau}^{}_{ra}(-\la_{ra}, k)\ts \eta,
\qquad \eta:=\theta_{\La'},
\end{equation}
where $\La'$ is the array obtained from $\La$ by increasing
each entry $\la_{aa},\dots,\la_{r-1,a}$ by $k$.  
We shall use a
(reverse) induction on the pairs $(r,a)$
ordered lexicographically, with the base
$\theta_{\La}=\zeta$. Note that
by the definition of admissible
arrays we must have $r\leq n-p+a$.

Identify the subalgebra of $\Y(n)$ generated by $t_{ij}(u)$
with $1\leq i,j\leq r$ with the Yangian $\Y(r)$.
By the induction hypothesis $\eta$ is a $\Y(r-1)$-singular vector
such that \eqref{tiieta1} holds with $n$ replaced by $r$, where
\beql{nunew}
\nu_i(u)=\begin{cases}
(u+\la_{ri})(u+\mu_i)\qquad&\text{if}\ \  i\leq a,\\
(u+\la_{r-1,i})(u+\mu_i)\qquad&\text{if}\ \  i> a.
\end{cases}
\non
\end{equation}
Therefore, if $\wt{\Tau}^{}_{ra}(-\la_{ra}, k)=\Tau^{}_{ra}(-\la_{ra}, k)$
then \eqref{biuthr} and \eqref{biu+1r} are immediate from
Lemma~\ref{lem:tii+1}.

Suppose now that 
$\wt{\Tau}^{}_{ra}(-\la_{ra}, k)=\Tau^{\ts\prime}_{ra}(-\la_{ra}, k)$.
Then $r=n-p+a$ and $k=k_a$. Therefore, $\la_{ra}=\la_a$ by \eqref{eqlala}.
It is clear from Lemma~\ref{lem:tii+1} that \eqref{biuthr} and \eqref{biu+1r}
hold for $i\ne a$ so that
we may assume $i=a$. 
In this case, due to \eqref{taanaze} and \eqref{taa+1naze},
it suffices to show that
\beql{taura=0}
\Tau_{ra}(-\la_{a},k_a)\ts\eta=0 
\end{equation}
and that for every $c=a+1,\dots,r$
the polynomial
\beql{polmu}
\nu_a(v)\ts\nu_{a+1}(v-1)\cdots\nu_{c-1}(v-c+a+1)\ts
\Tau^{}_{ra}(v+1,k_a-1)\ts\tau_{rc}(v-c+a)\ts\eta
\end{equation}
has zero of multiplicity at least two at $v=-\la_{a}$.

Suppose first that $p<n-1$ and
consider $\Tau_{ra}(-\la_{a},k_a)\ts\eta$.
By \eqref{qminor} we have
\beql{taunp}
\tau^{}_{ra}(v)=\sum_{\sigma\in\Sym_{r-a}}\sgn\sigma\cdot
T_{\sigma(a+1),a}(v)\cdots T_{\sigma(r),r-1}(v-r+a+1).
\end{equation}
By the induction hypothesis we have 
$T_{\sigma(r),r-1}(u)\ts\eta=0$
if $\sigma(r)<r-1$, while
\beql{tszeig}
T_{r-1,r-1}(u)\ts\eta
=(u+\la_{r-1,r-1})(u+\mu_{r-1})\ts\eta.
\non
\end{equation}
The factor $u+\mu_{r-1}$ is zero
if $u=-\la_{a}-r+a+1$ by \eqref{leqm}. Therefore, if $v=-\la_a$ then
we may assume that $\sigma(r)=r$ in \eqref{taunp} which gives
\beql{taulamus}
\tau^{}_{ra}(-\la_a)\ts\eta=\tau^{}_{r-1,a}(-\la_a)\ts T_{r,r-1}(-\mu_{r-1})\ts\eta.
\non
\end{equation}
Since $T_{r,r-1}(v)=\tau^{}_{r,r-1}(v)$ is a lowering operator,
we verify by an easy induction with the use of \eqref{center} and  
\eqref{taaa} that
\beql{simfinta}
\Tau^{}_{ra}(-\la_{a},k_a)\ts\eta=
\Tau^{}_{r-1,a}(-\la_{a},k_a)\ts \Tau^{}_{r,r-1}(-\mu_{r-1},k_a)\ts\eta.
\non
\end{equation}
We have $k_a=m_{r-1}-m_r$ by \eqref{leqm} and so,
the equality $\Tau_{ra}(-\la_{a},k_a)\ts\eta=0$ in the case $p<n-1$
will be implied by the fact that the vector
\beql{zeroyn2}
\wt{\eta}=T_{r,r-1}(-\mu_{r})\cdots T_{r,r-1}(-\mu_{r-1})\ts\eta
\end{equation}
is zero in $V(\la,\mu)$. 
Since the $\Y(n)$-module $V(\la,\mu)$ is irreducible, 
it will be sufficient to show, due to Proposition~\ref{prop:sing}, that
the vector $\wt{\eta}$
is annihilated by all operators	$B_i(u)$
with $i=1,\dots,n-1$. Since $T_{r,r-1}(u)$ commutes with
the lowering operators $\tau^{}_{ri}(v)$ we may assume that
the array $\La$ satisfies $\la_{ri}=\la_{r-1,i}$ for all $i\ne a$.
In other words, $\eta=\theta_{\La'}$ is a $\Y(r)$-singular vector
such that
\beql{tiiueta}
T_{ii}(u)\ts\eta=(u+\la_{ri})(u+\mu_i)\ts\eta\qquad\text{for}\quad
i=1,\dots,r.
\non
\end{equation}
By \eqref{tii+1na}
$\wt{\eta}$ is annihilated by the operators $T_{i,i+1}(u)$
for $i=1,\dots,r-2$, and by \eqref{tiinaze} and \eqref{taaa}
$\wt{\eta}$ is an eigenvector for the operators $T_{ii}(u)$
with $i=1,\dots,r-1$. We have $T_{r-1,r}(u)=[T_{r-1,r-1}(u),E_{r-1,r}]$.
By \eqref{ynknn-1}, $E_{r-1,r}\ts\wt{\eta}=0$
and therefore $T_{r-1,r}(u)\wt{\eta}=0$.
On the other hand, if $i\geq r$ then
$B_i(u)$ commutes with the elements
$T_{r,r-1}(v)$ by \eqref{center}.
Hence, by the induction hypothesis,
$B_i(u)\ts\wt{\eta}$ is a linear combination of the vectors
\beql{blinc}
T_{r,r-1}(-\mu_{r})\cdots T_{r,r-1}(-\mu_{r-1})\ts
\theta_{\La'+\delta_{ij}}
\non
\end{equation}
with $2\leq j\leq p$.
We conclude by induction on the weight of $\La'$ that
$B_i(u)\ts\wt{\eta}=0$ thus proving \eqref{taura=0}.

Consider now the polynomial \eqref{polmu}. Suppose first that $c=r$.
Note that $r-1>a$ since $p<n-1$.
We have
\beql{polnn}
\nu_{a}(v)
=(v+\la_{a})(v+\mu_{a}),\qquad
\nu_{r-1}(v)
=(v+\la_{r-1,r-1})(v+\mu_{r-1}).
\end{equation}
However, $\mu_{r-1}-r+a+1 =\la_{a}$ by \eqref{leqm}. 
This shows that the coefficient
of $\eta$ in \eqref{polmu} is divisible by $(v+\la_{a})^2$.
Let now $a+1\leq c< r$.	By \eqref{qminor},
\beql{taqqq}
\tau^{}_{rc}(-\la_{a}-c+a)=
\sum_{\sigma\in\Sym_{r-c}}\sgn\sigma\cdot
T_{\sigma(c+1),c}(-\la_{a}-c+a)\cdots T_{\sigma(r),r-1}(-\la_{a}-r+a+1).
\non
\end{equation}
We can repeat the argument which we have applied to the expression
\eqref{taunp} to show that the polynomial
$\Tau^{}_{ra}(v+1,k_a-1)\ts\tau_{rc}(v-c+a)\ts\eta$
has zero at $v=-\la_{a}$. Together with
\eqref{polnn} this completes the proof of \eqref{biuthr}
and \eqref{biu+1r} in the case $p<n-1$.

In the case $p=n-1$ we have $a=r-1$ and
\beql{zeroyn3}
\Tau_{r,r-1}(-\la_{r-1},k_{r-1})\ts\eta=
T_{r,r-1}(-\mu_{r})\cdots T_{r,r-1}(-\la_{r-1})\ts\eta.
\non
\end{equation}
Note that the operators $T_{r,r-1}(u)$ and $T_{r,r-1}(v)$ commute.
Therefore, due to \eqref{ineqml} it suffices to show that
the vector \eqref{zeroyn2} is zero. The argument
used in the case $p<n-1$ works here as well.

For $p=n-1$ the polynomial \eqref{polmu} equals
\beql{pomult2}
\nu_{r-1}(v)\ts\Tau_{r,r-1}(-\la_{r-1}+1,k_{r-1}+1)\ts\eta.
\end{equation}
If $\la_{r-1}=\mu_{r-1}$ then
$\nu_{r-1}(v)=(v+\la_{r-1})^2$. If $\la_{r-1}-\mu_{r-1}>0$ then
\beql{zeroyn4}
\Tau_{r,r-1}(-\la_{r-1}+1,k_{r-1}+1)\ts\eta=
T_{r,r-1}(-\mu_{r})\cdots T_{r,r-1}(-\la_{r-1}+1)\ts\eta.
\non
\end{equation}
Using again the fact that the vector \eqref{zeroyn2} is zero
we conclude that this vector is also zero.
In the both cases \eqref{pomult2} has zero of multiplicity at least two at
$v=-\la_{r-1}$ proving	\eqref{biuthr}
and \eqref{biu+1r}.

To prove \eqref{bathL2r} we note that
$B_i(u)$ commutes with $\Tau^{}_{sa}(v,k)$ for $i\geq s$ by \eqref{center}.
Therefore, it suffices to consider the case $i=r-1$. We derive from \eqref{qmrel} that
$B_{r-1}(u)=[A_{r-1}(u),E_{r-1,r}]$. Suppose first that $p<n-1$.
By \eqref{biuthr} and \eqref{biu+1r} the operator
$A_{r-1}(u)$ acts on $\theta_{\La}$
as multiplication by a polynomial in $u$. So it suffices to prove that
\beql{bathL2r2}
E_{r-1,r}\ts\theta_{\La}=
\sum_{j=1}^{r-1}\beta_{j}(\La)\ts\theta_{\La+\delta^{}_{r-1,j}},
\end{equation}
where $\beta_{j}(\La)$
are some constants.	Write $\theta_{\La}$ in the form \eqref{thlade}
and assume that $a<r-1$.
We now use Lemma~\ref{lem:ecomm}. By \eqref{ecomta} we have
\beql{er-1r}
E_{r-1,r}\ts\Tau^{}_{ra}(v,k)\ts\theta_{\La'}
=\Tau^{}_{ra}(v,k)\ts E_{r-1,r}\ts\theta_{\La'}
-k\ts\Tau^{}_{ra}(v,k-1)\ts
{T\ts}^{a+1\ts\cdots\ts r}_{a\ts\cdots\ts r-2,r}(v+k-1)\ts\theta_{\La'}.
\end{equation}
Further, \eqref{tautau}	gives
\beql{tautaur}
{T\ts}^{a+1\ts\cdots\ts r}_{a\ts\cdots\ts r-2,r}(u)\ts
{T\ts}^{a+1\ts\cdots\ts r-1}_{a+1\ts\cdots\ts r-1}(u)\ts\theta_{\La'}=
\tau^{}_{r-1,a}(u)\ts
{T\ts}^{a+1\ts\cdots\ts r}_{a+1\ts\cdots\ts r}(u)\ts\theta_{\La'}
+ \tau^{}_{ra}(u)\ts
{T\ts}^{a+1\ts\cdots\ts r-1}_{a+1\ts\cdots\ts r-2,r}(u)\ts\theta_{\La'}.
\end{equation}
By \eqref{qmrel}, 
${T\ts}^{a+1\ts\cdots\ts r-1}_{a+1\ts\cdots\ts r-2,r}(u)
=[{T\ts}^{a+1\ts\cdots\ts r-1}_{a+1\ts\cdots\ts r-1}(u),E_{r-1,r}]$.
Using the induction hypothesis we obtain
\beql{eihy}
E_{r-1,r}\ts\theta_{\La'}=
\sum_{j=a+1}^{r-1}\beta_{j}(\La')\ts\theta_{\La'+\delta^{}_{r-1,j}}.
\non
\end{equation}
On the other hand,
using \eqref{biuthr} and \eqref{biu+1r} we derive from \eqref{qminor} that 
\beql{tau+1r}
{T\ts}^{a+1\ts\cdots\ts r-1}_{a+1\ts\cdots\ts r-1}(u)\ts\theta_{\La'}
=\prod_{i=a+1}^{r-1}(u+l_{r-1,i}+a)(u+m_i+a)\ts \theta_{\La'},
\non
\end{equation}
while
\beql{tau+1rr}
{T\ts}^{a+1\ts\cdots\ts r}_{a+1\ts\cdots\ts r}(u)\ts\theta_{\La'}
=\prod_{i=a+1}^{r}(u+l_{ri}+a)(u+m_i+a)\ts \theta_{\La'},
\non
\end{equation}
since ${T\ts}^{a+1\ts\cdots\ts r}_{a+1\ts\cdots\ts r}(u)$
commutes with the lowering operators $\tau_{rb}(v)$. Thus, by
\eqref{tautaur}
\beql{taurfin}
\bal
{T\ts}^{a+1\ts\cdots\ts r}_{a\ts\cdots\ts r-2,r}(u)\ts\theta_{\La'}
{}=(u+l_{rr}+a)&(u+m_r+a)\prod_{i=a+1}^{r-1}\frac{u+l_{ri}+a}{u+l_{r-1,i}+a}
\ts\tau^{}_{r-1,a}(u)\ts\theta_{\La'}\\
{}+{}\sum_{j=a+1}^{r-1}{}&{}\frac{\beta_{j}(\La')}{u+l_{r-1,j}+a}\ts
\tau^{}_{ra}(u)\ts\theta_{\La'+\delta^{}_{r-1,j}}.
\eal
\end{equation}
Now put $v=-\la_{ra}$ and $k=\la_{ra}-\la_{r-1,a}$ into \eqref{er-1r}.
The denominator $u+l_{r-1,i}+a$ in \eqref{taurfin} becomes
$l_{r-1,i}-l_{r-1,a}$ at $u=v+k-1$.
Due to the conditions \eqref{leqm} and \eqref{eqlala}
the difference $l_{r-1,i}-l_{r-1,a}$ can only be zero if $i=a+1$.
Moreover, in this case $l_{r-1,a+1}=l_{r,a+1}=l_{a+1}$.	Then
$\La'+\delta^{}_{r-1,a+1}$ is not an admissible array so that
the summand with $j=a+1$ does not occur in the sum in \eqref{taurfin}.
The denominator $u+l_{r-1,i}+a$ with $i=a+1$ does not occur in the
product either, since it cancels	with $u+l_{r,a+1}+a$.
Thus the substitution $u=v+k-1$
into \eqref{taurfin}
is well defined. Using the fact that $\tau^{}_{r-1,b}(v)$ commutes
with $\tau^{}_{ra}(u)$ if $b\geq a$ we complete the proof of
\eqref{bathL2r2} for the case 
$\wt{\Tau}^{}_{ra}(v, k)=\Tau^{}_{ra}(v, k)$
in \eqref{thlade}.

Assume now that
$\wt{\Tau}^{}_{ra}(v, k)=
\Tau^{\ts\prime}_{ra}(v, k)$.
Then by \eqref{taitidef} we must have $r=n-p+a$ and $k=k_a=l_a-m_r$.
Moreover, we also have $\la_{ra}=\la_a$ by \eqref{eqlala}.
Take the derivative with respect to $v$ in \eqref{er-1r}
and put $v=-\la_a$. Note that the factor $u+m_r+a$ in \eqref{taurfin}
vanishes at $u=-\la_a+k_a-1$. Furthermore, as has been shown
above,
$\Tau^{}_{ra}(-\la_a, k_a)\ts\theta_{\La'+\delta^{}_{r-1,j}}=0$;
see \eqref{taura=0}. The application of the induction hypothesis
finally proves \eqref{bathL2r2} in the case $a<r-1$.

If $a=r-1$ in \eqref{thlade} then $\eta=\theta_{\La'}$ 
is a $\Y(r)$-singular vector, so that we may use
the relation \eqref{ynknn-1} 
to prove \eqref{bathL2r2}.
The same relation applies
in the case $p=n-1$,
where we also use the fact that the polynomial \eqref{pomult2} has zero 
of multiplicity at least two at
$v=-\la_{r-1}$.
\endproof

Consider the vector $\theta\in V(\la,\mu)$ defined in \eqref{singtau}.

\ble\label{lem:theta=0}
$\theta=0$.
\ele

\Proof
We shall be proving
by induction on the weight of $\La$ that
\beql{talze}
\Tau^{}_{n-p+1,1}(-\la_1,k_1)\ts\theta^{}_{\La}=0
\end{equation}
for all admissible arrays $\La$ such that the parameter
$r=r(\La)$ satisfies $r\geq n-p+2$.
For the induction base we note that
\beql{singze}
\Tau^{}_{n-p+1,1}(-\la_1,k_1)\ts\zeta=0.
\non
\end{equation}
Indeed, using \eqref{center} we find that the vector on the left hand side
is annihilated by the operators $B_i(u)$ with $i=n-p+1,\dots,n-1$.
On the other hand,
by \eqref{taaa+1}
it is also annihilated by 
$T_{i,i+1}(u)$ with $i=1,\dots,n-p-1$.
Further, $T_{n-p,n-p+1}(u)=[T_{n-p,n-p}(u), E_{n-p,n-p+1}]$
and we find from \eqref{tiinaze} and \eqref{ynknnne}
that it is also annihilated by $E_{n-p,n-p+1}$.
By Proposition~\ref{prop:sing} the vector must be zero.

Suppose now that $w(\La)\prec \la$. Denote the
left hand side of \eqref{talze} by $\wt{\theta}^{}_{\La}$.
We shall show that $B_i(u)\ts\wt{\theta}^{}_{\La}=0$ for all $i=1,\dots,n-1$.
By Lemma~\ref{lem:theLa}, the $\Y(n-p+1)$-span of
the vector $\theta^{}_{\La}$ is a highest weight module
with the highest weight defined from \eqref{biuthr} with $r=n-p+2$.
Exactly as above we find that $T_{n-p,n-p+1}(u)\ts\theta^{}_{\La}$ is zero.
Furthermore, the
operators $B_i(u)$ 
with $i=n-p+1,\dots,n-1$ commute with $\Tau^{}_{n-p+1,1}(-\la_1,k_1)$.
Therefore by \eqref{bathL2r}
$B_i(u)\ts\wt{\theta}^{}_{\La}$
is a linear combination of the vectors 
$\Tau^{}_{n-p+1,1}(-\la_1,k_1)\ts\theta^{}_{\La+\delta_{ij}}$.
If $i\geq n-p+2$ then the arrays $\La+\delta_{ij}$ satisfy
the condition $r\geq n-p+2$ on the parameter $r=r(\La)$ used in Lemma~\ref{lem:theLa}, 
and we complete the proof in this case
applying the induction hypothesis.

If $i=n-p+1$ then
\beql{septs} 
\theta^{}_{\La+\delta_{n-p+1,j}}=\tau^{}_{n-p+1,j}(-\la_{n-p+1,j}-1)\ts
\theta^{}_{\La'},
\non
\end{equation}
for some array $\La'$ for which the corresponding parameter $r'=r(\La')$
satisfies $r'\geq n-p+2$.
However, $\tau^{}_{n-p+1,j}(v)$ is permutable with $\Tau^{}_{n-p+1,1}(-\la_1,k_1)$
which again ensures that $B_i(u)\ts\wt{\theta}^{}_{\La}=0$
by the induction hypothesis.
\endproof

By \eqref{copr} and \eqref{epi0} the operators
$T_{ij}(u)=u^2\ts t_{ij}(u)$ in $L(\la)\ot L(\mu)$ are polynomials in $u$.
Therefore, we can
introduce the vector $\wt{\theta}\in L(\la)\ot L(\mu)$ by
\beql{thetot}
\wt{\theta}=\Tau^{}_{n-p+1,1}(-\la_1,k_1)\ts
\Tau^{\ts\prime}_{n-p+2,2}(-\la_2,k_2)\cdots 
\Tau^{\ts\prime}_{np}(-\la_p,k_p)\ts(\xi\ot\xi'),
\non
\end{equation}
cf. \eqref{singtau}.
Here $\xi$ and $\xi'$ are the highest vectors of
the $\gl_n$-modules $L(\la)$ and $L(\mu)$, respectively.

\ble\label{lem:n-1ot} $\wt{\theta}\ne 0$.
\ele

\Proof Write the vector $\wt{\theta}$ in the form
\beql{gtLM}
\wt{\theta}=\sum_{\La,M}c^{}_{\La,M}\ts\xiL\ot\xi'_M,
\end{equation}
where $\xiL$ and $\xi'_M$ are the Gelfand--Tsetlin basis vectors in $L(\la)$ and
$L(\mu)$. It suffices to show that at least one coefficient 
$c^{}_{\La,M}$ is nonzero.
We shall calculate these coefficients for the case 
$\xi'_M=\xi'_{M^0}$ is the highest
vector $\xi'$ of $L(\mu)$. That is, the $(k,i)$-th entry of the pattern $M^0$
coincides with $\mu_i$ for all $k$ and $i$. 
We prove by a (reverse) induction on $a$ that for any $a=2,\dots,p$ we have
\beql{thethi}
\Tau^{\ts\prime}_{n-p+a,a}(-\la_a,k_a)\cdots 
\Tau^{\ts\prime}_{np}(-\la_p,k_p)\ts(\xi\ot\xi')
=\sum_{\La,M}c^{(a)}_{\La,M}\ts\xiL\ot\xi'_M,
\end{equation}
where for each $M$ occurring
in this expansion,
\beql{comu}
\mu-w(M)=\sum_{i=a}^{n-1} q_i\ts (\ve_i-\ve_{i+1}),\qquad q_i\in\ZZ_+,
\end{equation}
and 
$c^{(a)}_{\La^{(a)},M^0}\ne 0$ for the pattern $\La^{(a)}$
defined for each $a=1,\dots,p$ as follows. 
The entry $\la^{(a)}_{si}$ of $\La^{(a)}$
coincides with $\la_i$ unless $i=a,\dots,p$ and $s<n-p+i$.
For the entries with these $s$ and $i$ we have $\la^{(a)}_{si}-i+1=m_{n-p+i}$.
The betweenness conditions \eqref{betw}
for $\La^{(a)}$ are guaranteed by the assumptions
\eqref{ineqml} and \eqref{leqm}.
 
Suppose that $a\leq p$ and denote the left hand side of
\eqref{thethi} by $\theta^{(a)}$.
By Proposition~\ref{prop:delqm}, we have
\beql{deltau}
\Delta(\tau_{ra}(v))=\sum_{b_1<\cdots<b_{n-p}}
{T\ts}^{a+1\ts\cdots\ts  r}_{b_1\ts\cdots\ts b_{n-p}}(v)
\ot {T\ts}^{b_1\ts\cdots\ts b_{n-p}}_{a\ts\cdots\ts  r-1}(v),
\qquad r:=n-p+a,
\non
\end{equation}
where the quantum minor operators in $L(\la)$ and $L(\mu)$ are defined
by the formulas \eqref{qminor} where the $t_{ij}(u)$ are replaced
with the polynomial operators $T_{ij}(u)=u\ts t_{ij}(u)$.
If $w$ is a weight of the $\gl_n$-module $L(\mu)$ then
$w\preceq \mu$. Therefore, by the induction hypothesis,
if $b_1<a$ then
\beql{theio}
{T\ts}^{b_1\ts\cdots\ts b_{n-p}}_{a\ts\cdots\ts  r-1}(v)\ts \xi'_M=0
\non
\end{equation}
for each pattern $M$ occurring in the expansion \eqref{thethi} for
$\theta^{(a+1)}$.
This proves \eqref{comu}. Furthermore, for any $k\in\ZZ_+$
the tensor products of the form $\xiL\ot\xi'_{M^0}$ which occur in
the expansion of $\Tau^{}_{ra}(v,k)\ts\theta^{(a+1)}$ 
should have the form
\beql{compo}
\Tau^{}_{ra}(v,k)\ts\xi^{}_{\La^{(a+1)}}\ot 
{T\ts}^{a\ts\cdots\ts  r-1}_{a\ts\cdots\ts  r-1}(v+k-1)\cdots
{T\ts}^{a\ts\cdots\ts  r-1}_{a\ts\cdots\ts  r-1}(v)\ts\xi'.
\end{equation}
The coefficient of $\xi'$ equals
\beql{ceffxip}
\prod_{j=1}^k (v+\mu_a+j-1)(v+\mu_{a+1}+j-2)\cdots (v+\mu_{r-1}+j-n+p).
\non
\end{equation}
The conditions \eqref{ineqml} and \eqref{leqm} imply that for 
$k=k_a$ there is a unique factor in this product which vanishes at 
$v=-\la_a$.
Therefore, the derivative of
\eqref{compo} with $k=k_a$ at $v=-\la_a$ is, up to a nonzero constant factor,
\beql{dera}
\Tau^{}_{n-p+a,a}(-\la_a,k_a)\ts\xi^{}_{\La^{(a+1)}}\ot\xi'.
\non
\end{equation}
By the definition \eqref{xildef}, this coincides with
$\xi^{}_{\La^{(a)}}\ot\xi'$ which proves \eqref{thethi}.

Similarly, the application of the operator 
$\Tau^{}_{n-p+1,1}(-\la_1,k_1)$ to the vector $\theta^{(2)}$
produces a linear combination \eqref{gtLM}. Here
the coefficient $c^{}_{\La^{(1)},M^0}$
is the product of $c^{(2)}_{\La^{(2)},M^0}$
and the factor
\beql{ceffxi1}
\prod_{j=1}^k (-\la_1+\mu_1+j-1)(-\la_1+\mu_2+j-2)\cdots (-\la_1+\mu_{n-p}+j-n+p)
\non
\end{equation}
with $k=l_1-m_{n-p+1}$ which comes from the expansion of the
coefficient of $\xi'$ in \eqref{compo} for $a=1$.
It remains to note that by \eqref{lone} (which holds for all $p\leq n-1$)
this factor is nonzero. \endproof

If the $\Y(n)$-module $L(\la)\ot L(\mu)$ is irreducible then by
\eqref{copr} and \eqref{epi0} it is isomorphic to the
highest weight module $V(\la,\mu)$.
Lemmas~\ref{lem:theta=0} and \ref{lem:n-1ot} therefore imply that
this contradicts to the assumption
\eqref{pcond}, thus proving Theorem~\ref{thm:necc}.

Theorem~\ref{thm:main1} is now implied by Theorems~\ref{thm:sufc}
and \ref{thm:necc} due to Proposition~\ref{prop:equiv}.


\end{document}